\documentclass{article}

\usepackage{amsmath}
\usepackage{amsfonts}
\usepackage{amscd}
\usepackage{amssymb}

\input xypic

\newtheorem{defn}{Definition}[section]
\newtheorem{thm}[defn]{Theorem}
\newtheorem{prop}[defn]{Proposition}
\newtheorem{lemma}[defn]{Lemma}
\newtheorem{eg}[defn]{Example}

\newtheorem{cor}[defn]{Corollary}

\newcommand{\lm}{\ensuremath{\longrightarrow}}

\newcommand{\eps}{\varepsilon}

\DeclareMathOperator{\Rhom}{\mbox{RHom}}
\DeclareMathOperator{\Hom}{\mbox{Hom}}
\DeclareMathOperator{\shom}{\ensuremath{\mathcal{H}\mathit{om}}}

\DeclareMathOperator{\send}{\ensuremath{\mathcal{E}\!\mathit{nd}}}
\DeclareMathOperator{\End}{\mbox{End}}
\DeclareMathOperator{\Aut}{\mbox{Aut}}
\DeclareMathOperator{\Ext}{\mbox{Ext}}

\DeclareMathOperator{\des}{\mbox{\footnotesize res}}
\DeclareMathOperator{\spec}{\mbox{Spec}\,}

\DeclareMathOperator{\mo}{\mbox{mod}}

\DeclareMathOperator{\Br}{\mbox{Br}\,}
\DeclareMathOperator{\Pic}{\mbox{Pic}\,}

\DeclareMathOperator{\Hilb}{\mbox{Hilb}\,}

\DeclareMathOperator{\s}{\sigma}
\DeclareMathOperator{\e}{\varepsilon}

\DeclareMathOperator{\w}{\omega}
\DeclareMathOperator{\z}{\zeta}

\DeclareMathOperator{\lan}{\langle}
\DeclareMathOperator{\ran}{\rangle}

\DeclareMathOperator{\Q}{\mathbb{Q}}
\DeclareMathOperator{\Z}{\mathbb{Z}}

\DeclareMathOperator{\calo}{\mathcal{O}}

\DeclareMathOperator{\ow}{\mathcal{O}_{W}}

\DeclareMathOperator{\oy}{\mathcal{O}_{Y}}

\DeclareMathOperator{\At}{\ensuremath{\tilde{A}}}

\DeclareMathOperator{\Wt}{\ensuremath{\tilde{W}}}
\DeclareMathOperator{\Xt}{\ensuremath{\tilde{X}}}
\DeclareMathOperator{\Yt}{\ensuremath{\tilde{Y}}}
\DeclareMathOperator{\Zt}{\ensuremath{\tilde{Z}}}
\DeclareMathOperator{\Gbar}{\ensuremath{\overline{G}}}

\begin{document}

\begin{center}
  \LARGE \textbf{McKay Correspondence for Canonical Orders}
\end{center}

\begin{center}
  DANIEL CHAN\footnote{Supported by ARC Discovery project grant DP0557228}
\end{center}

\begin{center}
  {\em University of New South Wales}
\end{center}

\begin{center}
e-mail address:{\em danielc@unsw.edu.au}
\end{center}

\begin{abstract}
Canonical orders, introduced in the minimal model program for orders [CI05], are simultaneous generalisations of Kleinian singularities $k[[s,t]]^G, G < SL_2$ and their associated skew group rings $k[[s,t]]*G$. In this paper, we construct minimal resolutions of canonical orders via non-commutative cyclic covers and skew group rings. This allows us to exhibit a derived equivalence between minimal resolutions of canonical orders and the skew group ring form of the canonical order in all but one case. The Fourier-Mukai transform used to construct this equivalence allows us to make explicit, the numerical version of the McKay correspondence for canonical orders noted in [CHI], which relates the exceptional curves of the minimal resolution to the indecomposable reflexive modules of the canonical order. 
\end{abstract}

We assume throughout this paper that everything is over an algebraically closed base field $k$ of characteristic zero.

\section{Introduction}  \label{sintro} 

There has been a proliferation of research on the McKay correspondence over the last couple of decades. This beautiful theory links the theory of singularities in algebraic geometry with non-commutative algebra and group theory. There are now many incarnations and generalisations of McKay's original correspondence, and we present yet another one here based on the concept of canonical orders introduced in [CHI]. 

Let us review one of the classical versions of the McKay correspondence as it would be most relevant for us. Consider a finite subgroup $G < SL_2$ which has a natural action on the power series ring $k[[s,t]]$. Then $\spec k[[s,t]]^G$ is a Kleinian (or canonical surface) singularity and the McKay correspondence states that there is a bijection between the set of exceptional curves in a minimal resolution for $\spec A$ and the set of indecomposable reflexive $A$-modules not isomorphic to $A$. It is customary now, following Kapranov-Vasserot [KV] (and presumably Verdier), to view this ``numerical McKay correspondence'' as a consequence of a derived equivalence as follows. Firstly, note that the order $A:= k[[s,t]]*G$ is reflexive Morita equivalent to $k[[s,t]]^G$ so, in particular, has the same number of indecomposable reflexive modules. Van den Bergh [V] and others interpret $A$ as a non-commutative resolution of $\spec k[[s,t]]^G$ since it has global dimension two (which suggests smooth) and the order $A$ has centre $k[[s,t]]^G$ and is Azumaya away from the closed point (which suggests birational to $k[[s,t]]^G$). Kapranov-Vasserot show there is a derived equivalence between $A$ and the minimal resolution of $\spec k[[s,t]]^G$. 

In [CI05], a version of Mori's minimal model program was introduced for orders on surfaces. As in the Mori program, there are non-commutative analogues for canonical singularities dubbed, canonical orders (the definition of which is reviewed in \S~\ref{sram}). They include as examples, Kleinian singularities $k[[s,t]]^G, G < SL_2$ as well as their associated skew group rings $k[[s,t]]*G$. Unlike other non-commutative generalisations of Kleinian singularities, there is a notion of a minimal resolution for a canonical order $B$. This minimal resolution is a certain order on a resolution of the centre of $B$ (see section~\ref{sskewres} for the definition). Naturally, we used the minimal resolution in [CHI] to study these orders. We showed in particular, that any canonical order is reflexive Morita equivalent to another canonical order of the form $A:=\eps k[[s,t]]*G'$ where $G'$ is a finite group and $\eps$ is a central primitive idempotent of $k[[s,t]]*G'$. We interpret this order $A$ as a non-commutative resolution of the original canonical order. It is thus natural to ask if $A$ is derived equivalent to the minimal resolution of the original order. In this paper, we show that this is true in all except possibly one case. 

\begin{thm}   \label{tMcKay} 
Given a canonical order $B$ of any ramification type other than $DL_n$ (see definition~\ref{dcanonical}), there is a reflexive Morita equivalent order of the form $A=\eps k[[s,t]]*G'$ such that $A$ is derived equivalent to the minimal resolution of $B$.
\end{thm}

The key to proving this theorem lies in studying the minimal resolution of a canonical order. The ramification data of the resolution is well-known from [CHI], but up till now, nothing else was known about the resolution. We remedy this by showing that terminal resolutions can always be constructed as non-commutative cyclic covers (as introduced in [C]). In the process, we see that several terminal resolutions of canonical orders can also be realised as skew group rings. This allows one to use an  equivariant version of the classical ``commutative'' McKay correspondence to prove the above theorem in all cases except those of ramification type $L_n$ and $DL_n$. To prove the theorem in the case $L_n$, we use the description of the minimal resolution as a non-commutative cyclic cover and construct explicitly, the kernel of the Fourier-Mukai transform which exhibits the desired derived equivalence. We unfortunately, do not know how to handle the case $DL_n$. 

The ``derived'' McKay correspondence we use is different from the one in [VdB] since the minimal resolution of an order is non-commutative. When the minimal resolution is a skew group ring, it is almost a special case of Kawamata's equivariant derived McKay correspondence [Kaw]. Hence, the content of the results here, is that canonical orders provide interesting examples of derived equivalences. The derived equivalence for type $L_n$ seems to be quite different from other known examples. 

The next natural question to ask is whether or not one can extract a numerical form of the McKay correspondence for canonical orders using theorem~\ref{tMcKay}. In [CHI] already, such a numerical McKay correspondence was observed and can be roughly stated as follows. If $B$ denotes a canonical order, the set of indecomposable reflexive $B$-modules can be partitioned into certain orbits containing one or two modules (see definition of permissible modules in [CHI, section~8]). Then [CHI, theorem~8.6] the number of these orbits is one more than the number of exceptional curves in a minimal resolution of $B$. In this paper, we explain this numerical McKay correspondence for types other than $DL_n$ by using the derived equivalence of theorem~\ref{tMcKay}. As a result, we will actually exhibit a bijection between the exceptional curves and the orbits of indecomposable reflexives not containing $B$. Furthermore, it allows us to give a more enlightening version of the numerical McKay correspondence in terms of the ramification type of the exceptional curves (see the paragraph preceding proposition~\ref{pskewgroup} for definitions). 

\begin{thm}  \label{tnumMcKay} 
Let $E_1,\ldots,E_r$ be the exceptional curves of the minimal resolution of a canonical order $B$. The number of indecomposable reflexive $B$-modules up to isomorphism is $n_0 + n_1 + \ldots + n_r$ where 
\begin{enumerate}
\item  $n_0 = 2$  if $B$ is of type $A_{1,2,\z},BL_n$ or $B_n$ and is 1 otherwise (see proposition~\ref{ptrivialcase} for a more uniform explanation).
\item $n_i=1$ if $E_i$ is of type $0,I_{1,e},C_2$ or $X_2$.
\item $n_i=2$ if $E_i$ is of type $I_{2,e}$. 
\end{enumerate}
\end{thm}

In section~\ref{sram} we remind the reader of the definition of terminal and canonical orders in terms of ramification data. To simplify the treatment and reduce the reliance on [CHI], we have opted for a non-standard definition of canonical orders. Section~\ref{sskew} is devoted to constructing terminal and canonical orders as skew group rings in the complete local case. We define minimal resolutions of canonical orders in section~\ref{sskewres} and examine when they can also be constructed via skew group rings. The next two sections describe an equivariant version of the derived McKay correspondence and a procedure for extracting the numerical McKay correspondence in the cases of interests. Section~\ref{sgeneral} shows how terminal resolutions can all be described, up to Morita equivalence, as non-commutative cyclic covers. The last section forms a  case-by-case analysis depending on the ramification type of the canonical order. We construct for each type, the minimal resolution as a non-commutative cyclic cover and, where possible, as a skew group ring. In the cases where there is a description via a skew group ring, we apply the results concerning the equivariant derived McKay correspondence in sections~\ref{sequi} and \ref{snumMcKay} to prove theorem~\ref{tMcKay} and elucidate the numerical McKay correspondence as discovered in [CHI].

\textbf{Acknowledgements:} I would like to thank Colin Ingalls and Paul Hacking for helpful discussions. In particular, I learnt about Kawamata's equivariant derived equivalence from Paul Hacking.

\section{Ramification of terminal and canonical orders}  \label{sram} 

Terminal and canonical orders are certain orders on surfaces which arise naturally in the minimal model program for orders over surfaces as studied in [CI05]. They are defined essentially in terms of certain geometric invariants called ramification data. In this section, we review the ramification data of terminal and canonical orders. Full details of the theory can be found in [CI05] and [CHI], and we revise here but the bare minimum.

We study the complete local case first. Let $R$ be a noetherian complete local normal domain of dimension two. An $R$-order $A$ in a central simple $k(R)$-algebra is {\em normal} [CI05, definition~2.3] if it is reflexive as an $R$-module and for every prime divisor $C \subset \spec R$ we have that 
\begin{enumerate}
\item the localisation at $C$, $A_C$ is hereditary.
\item the Jacobson radical $\mbox{rad}\, A_C$ is a principal left and right ideal.
\end{enumerate}

We note here that maximal orders are normal and skew group rings tend to be normal too. The condition is also stable under \'etale base change. The importance of normality for us is the following result.

\begin{prop}  \label{pnormalram} 
Given a normal $R$-order as above and a prime divisor $C \subset \spec R$, then $\tilde{k}(C):= Z(A_C/\mbox{rad}\, A_C)$ is a cyclic extension of $k(C)$ of the form $L^n$ where $L$ is a cyclic field extension of $k(C)$.
\end{prop}
Above, $Z$ denotes the centre of a ring as usual. The {\em ramification index} of $A$ at $C$ is 
\[   e_C:= \dim_{k(C)} \tilde{k}(C)                  . \]
This gives us the primary ramification of the order. We say $A$ {\em ramifies} at $C$ if $e_C > 1$ and say the corresponding curve is a {\em ramification curve}. The {\em ramification data} are the centre $Z$, the set of ramification curves $\{C\}$ and the corresponding field extensions $\tilde{k}(C)/k(C)$. 

Note that the extension $\tilde{k}(C)/k(C)$ corresponds to a cyclic cover $\tilde{D}_C/D_C$ of smooth (but possibly reducible) curves. In the complete local case, this cover is determined by its ramification and, following [AdJ], we will refer the ramification of $\tilde{D}_C/D_C$ to be secondary ramification of the order $A$. Thus the ramification data reduces to the primary and secondary ramification data. 

\begin{defn} \label{dterminal} 
([CI05, definition~2.5]) A normal $R$-order $A$ is {\em terminal} if its ramification data satisfies the following conditions.
\begin{enumerate}
\item $R$ is smooth.
\item the ramification curves are smooth and cross normally (if at all).
\item secondary ramification occurs only where ramification curves intersect.
\item if there are two ramification curves $C_1,C_2$ with ramification indices $e_1,e_2$ then, (swapping indices if necessary) we have $e_1|e_2$ and the secondary ramification of both curves is $e_1$ at their point of intersection.
\end{enumerate}
\end{defn}
In other words, $A$ is terminal if we can find coordinates $u,v \in R$ such that $R \simeq k[[u,v]]$ where furthermore i) $A$ is unramified away from the normal crossing lines $uv=0$, ii) there are integers $e,n \geq 1$ such that the ramification indices at $u=0$ and $v=0$ are $e$ and $ne$ respectively, and finally iii) the secondary ramification indices of both lines is $e$ at $(u,v)=(0,0)$.

\begin{defn}  \label{dcanonical} 
A normal $R$-order $A$ is {\em canonical} if it has one of the ramification data in the table below. 
\end{defn}

\begin{tabular}{c|c|c|c|c}
type  & $R$  & ramification curves  & $e_C$  & $e_p$ \\ \hline
$A_{1,2,\xi}$  & $k[[u,v]]$  & $u=0$      & $2e$      & $e$     \\ 
              &             & $v=0$      & $2e$      & $e$     \\\hline
$BL_n$  & $k[[u,v]]$  & $v^2=u^{2n+1}$  &  2  &  1      \\\hline
$B_n$  & $k[[u,v]]$  & $v=u^n$      & 2      & 1     \\
              &             & $v=-u^n$      & 2      & 1     \\\hline
$L_n$  & $k[[u,v]]$  & $v=u^{n+1}$      & 2      & 2     \\
              &             & $v=-u^{n+1}$      & 2      & 2     \\\hline
$DL_n$  & $k[[u,v]]$  & $u=0$      & 2      & 2     \\
              &             & $v^2=u^{2n-1}$      & 2      & 2     \\\hline
$BD_n$  & $k[[u,v]]$  & $u=0$      & 2      & 2     \\
              &             & $v=u^{n-1}$      & 2      & 2     \\
              &             & $v=-u^{n-1}$      & 2      & 1     \\\hline
$A_{n,\xi}$  & $k[[u,v,w]]/(uv-w^{n+1})$   & $w=0=u$   & $e$   & $e$   \\
             &                            & $w=0=v$   & $e$   & $e$   \\\hline
$A^1_n,D^1_n,E^1_{6,7,8}$  & klein sing $k[[x,y]]^G, G < SL_2$  & $\varnothing$    &  & \\ \hline
\end{tabular}

\vspace{2mm}
In the table above, $R$ is the centre of the order, the ramification curves are listed for each type with associated primary ramification indices under $e_C$ and secondary ramification at the closed point $p \in \spec R$ listed under $e_p$. The last line of the table represents orders on a Kleinian singularity which are unramified (in codimension one). The original definition of canonical orders is given in [CI05; page~444] and is based on the concept of discrepancy for orders. That it is equivalent to the definition above is a theorem in [CHI; theorem~6.5].

\section{Skew group rings: local case}   \label{sskew} 

In this section and the next, we address the question: given terminal or canonical ramification data, how do you construct a skew group ring with that ramification data? We give here a recipe we learnt from work of Artin [A86, section~5] and Le Bruyn-Van den Bergh-Van Oystaeyen [LVV, section~5], though presumably it was known in some form or another beforehand. More details are also also reviewed in [CHI, section~4]. The recipe works effectively when the centre is a complete local ring. The construction in this case, is given in this section. The next section will deal with constructing minimal resolutions of canonical orders via skew group rings. 

Before reviewing the recipe, we note that in the cases of interest, the Brauer group of the centre is trivial, so terminal orders are determined up to Morita equivalence by their ramification data (see theorem~\ref{tcyclicexist} for a precise statement as well as [CI05, theorem~2.12 and corollary~2.13]). Hence, from the point of view of non-commutative algebraic geometry, where geometric objects manifest themselves via their category of quasi-coherent sheaves, constructing a terminal order amounts to constructing an order with the same ramification data. Unfortunately, a similar uniqueness result does not hold for canonical orders. Nevertheless, in the complete local case we know by [CHI, lemma~3.2], that they are all reflexive Morita equivalent to a (direct factor of a) skew group ring $\eps k[[s,t]]*G'$ ($\eps$ a primitive central idempotent) where $G'$ is a group determined uniquely by ramification data. Furthermore, [CHI, lemma~3.2] also shows this skew group ring has the same ramification data as the original canonical order and may be thought of as a ``non-commutative resolution'' of the order.

The key to the recipe lies in finding a smooth Galois cover of the surface with the same ramification as the primary ramification data. To be more precise, let us restrict our attention for the rest of this section to the complete local case and let $R$ be a normal noetherian complete local domain of dimension two. Consider the ramification data $\{R,e_C,e_p\}$ of a terminal or canonical order as listed in section~\ref{sram}. By [CHI, proposition~5.1], there is a Galois cover $k[[s,t]]/R$ whose ramification indices coincide with the $e_C$'s. Let $G$ be the Galois group of $k[[s,t]]/R$.

\begin{eg}  \label{eartincover}
Case of terminal ramification.
\end{eg}
Consider the ramification data $R=k[[u,v]]$, the two ramification curves are $u=0,v=0$ with respective ramification indices $e,ne$ and secondary ramification $e$ for both these curves at $(u,v)=(0,0)$. Up to isomorphism, terminal orders have such ramification data. Let $\zeta \in k$ be a primitive $ne$-th root of unity. In this case, we let 
\[ G:= \langle \bar{\s},\bar{\tau} | \bar{\s}^e=1=\bar{\tau}^{ne}, \bar{\s}\bar{\tau}=\bar{\tau}\bar{\s} \rangle 
  \simeq \Z / e\Z \times \Z / ne \Z \]
act on $k[[s,t]]$ by $\bar{\s}: s \mapsto \zeta^n s, t \mapsto t$ and $\bar{\tau}:s \mapsto s, t \mapsto \zeta t$. We note that $k[[s,t]]^G = k[[u,v]]$ if we identify $u=s^e,v=t^{ne}$. Hence, $k[[s,t]]/R$ is the Galois cover with the above ramification data. 

\vspace{2mm}

Returning to the general setup, we need a way to incorporate the secondary ramification data. Consider now a central extension of $G$ by the group $\mu_e$ of $e$-th roots of unity. 
\[   1 \lm \mu_e \lm G' \lm G \lm 1    \]
Such extensions are classified by $H^2(G,\mu_e)$. In the applications, $e$ will be the largest value of $e_p$. We can extend the $G$-action on $k[[s,t]]$ to $G'$ by letting $\mu_e$ act trivially. This means that the skew group ring $k[[s,t]]*G'$ contains a central subalgebra $k\mu_e$. Let $\rho$ be a generator for $\mu_e$. Now $k\mu_e \simeq k\times k \times \ldots \times k$, the product of $e$ copies of the field $k$. It contains $e$ primitive idempotents. For example 
\begin{equation} \label{eidemp} %\marginpar{eidemp}
 \varepsilon = \frac{1}{e}(\rho^{e-1} + \zeta^n\rho^{e-2} + \ldots + \zeta^{n(e-1)}) 
\end{equation} 
is the idempotent corresponding to the factor where $\rho$ acts as the primitive $e$-th root of unity $\zeta^n$. Now $\eps$ is also a central idempotent of $k[[s,t]]*G'$ so we may consider the direct factor 
\[\eps k[[s,t]]*G' \simeq \frac{k[[s,t]]*G'}{(\rho-\zeta^n)}  \]
which must also have global dimension two. This algebra can also be written as a cross product algebra. From this viewpoint, we thus see that if $\alpha \in H^2(G,\mu_e)$ corresponds to the extension $G'$ then, $\eps k[[s,t]]*G'$ depends only on the image $\beta$ of $\alpha$ in $H^2(G,k^*)$. Now proposition~4.9 of [A86] shows that $\beta \in H^2(G,k^*)$ determines the Brauer class of the central simple algebra $\eps k((s,t))*G'$ and hence, the secondary ramification of $\eps k[[s,t]]*G'$. To simplify terminology, we shall refer to direct factors of skew group rings as skew groups rings also.

\begin{eg}   \label{ecentralext}  
Case of terminal ramification.
\end{eg}
We continue the notation in example~\ref{eartincover}. We compute easily enough that $H^2(G,k^*)\simeq \mu_e$ for $G = \Z / e\Z \times \Z / ne \Z $. We let $G'$ be a central extension corresponding to a generator of the cyclic group $H^2(G,k^*)$. In fact, we may assume
\[ G' = \langle \s,\tau,\rho | \s^e=1=\tau^{ne} = \rho^e, \s\tau=\rho\tau\s,\rho\s=\s\rho,
\rho\tau=\tau\rho \rangle \]
and the natural map $G'\lm G$ sends $\s \mapsto \bar{\s}, \tau \mapsto \bar{\tau}, \rho \mapsto 1$. 

\begin{prop}   \label{pskewterminal}  
Let $G'$ be the central extension of $\Z / e\Z \times \Z / ne \Z$ by $\mu_e$ defined above. Let $G'$ act on $k[[s,t]]$ as in example~\ref{eartincover}. Then $\eps k[[s,t]]*G' \simeq k[[s,t]]*G'/(\rho-\zeta^n)$ is a terminal order with ramification data:
\begin{enumerate}
\item the centre is $R=k[[u,v]]$ where $u=s^e,v=t^{ne}$.
\item the ramification curves are $u=0,v=0$ with corresponding ramification indices $e_C=e,ne$.
\item the secondary ramification of both curves is $e$ at the node $u=0=v$.
\end{enumerate}
\end{prop}
\textbf{Proof.} The theory of crossed product algebras shows that $A=\eps k[[s,t]]*G'$ is an order with centre  $k[[u,v]]$ which is unramified away from $uv=0$. We check the ramification at $u=0$ first which involves computing the residue ring extension $A_P/ \mbox{rad}\, A_P$ where $P=(u)$. Now $s$ is normal and nilpotent in $A/(u)$ so lies in $\mbox{rad}\, A_P$. Furthermore, $A/(s) \simeq \eps k[[t]]*G'$ so the fact that $\eps k((t))*G'$ is semisimple ensures 
\[   A_P/\mbox{rad}\, A_P = \eps k((t))*G'   .\]
which has centre $k((t^{ne}))[t^n\s^{-1}]$ which is totally ramified over $k(P)=k((t^{ne}))$ with ramification index $e$ as was to be shown. Observe also that $\mbox{rad}\, A_P$ is principal generated byt the normal element $s$.

Similarly, at $P=(v)$ we find 
\[   A_P/ \mbox{rad}\, A_P = \eps k((s))*G'   .\]
which has centre $k((s))[s\tau^{-1},\tau^e] \simeq k((s))[s\tau^{-1}]^n$. This has degree $ne$ over $k(P)=k((s))$ and the secondary ramification index is $e$. Again we see that $\mbox{rad}\, A_P$ is principal. Since $A$ is Azumaya away from $uv=0$, we see that $A$ is also a normal order. 

\vspace{2mm} 

The above proposition constructs for each ramification type of a terminal order, a skew group ring with that ramification. The following well-known proposition is useful for studying skew group rings. 

\begin{prop}   \label{pdecomp} 
Let $Z$ be an affine scheme and $G$ be a finite group acting faithfully on $\calo_Z$ and transitively on the connected components of $Z$. Let $Y \subset Z$ be a connected component and $D \leq G$ be the corresponding decomposition group. Let $G'$ be a central extension of $G$ by $\mu_e$,  $D'$ the corresponding extension of $D$ and $\eps \in k\mu_e$ be a primitive idempotent. Then 
\[\eps \calo_Z * G' \simeq (\eps \oy*D')^{n \times n}  \]
where $n = [G:D]$. 
\end{prop}
\textbf{Proof.} Since $G$ acts transitively on connected components, $\calo_Z$ is a direct product of copies of $\oy$. If $\{\eps_i\}$ are the characteristic functions of the connected components, then there is a Peirce decomposition of $\eps \calo_Z *G'$ corresponding to the idempotents $\{\eps_i\}$. We omit the elementary verification that this gives the isomorphism stated in the proposition.
\vspace{2mm} 

%\begin{prop}   \label{pskewsplitterm}  \marginpar{pskewsplitterm}
%Let $\tau$ be an involution on a commutative ring $R$. We let $\tau$ act on $R \times R$ diagonally and %$\s$ be the automorphism of $R \times R$ which maps $(r,r') \mapsto (r',r)$. Suppose $G'$ is the central extension of $G=\langle \s \rangle \times \langle \tau \rangle$ in example~\ref{ecentralext}. Then 
%\[  \frac{(R \times R)*G'}{(\rho+1)} \simeq (R*\langle \tau \rangle)^{2 \times 2}   .\]
%\end{prop}
%\marginpar{need to generalise}
%\textbf{Proof.} The isomorphism is given by mapping $R \times R$ to the corresponding diagonal matrix and 
%\[  \tau \mapsto 
%\begin{pmatrix}
% \tau & \\ & -\tau
%\end{pmatrix}, \hspace{5mm}
%\s \mapsto 
%\begin{pmatrix}
% & 1 \\ 1 & 
%\end{pmatrix}
%.\]

The above recipe was used in [CHI] to construct a skew group ring $\eps k[[s,t]]*G'$ for each of the ramification types of a canonical order. We merely list the answer here by specifying $G'$ and its action on $k[[s,t]]$. It turns out that in almost all cases, the extension $G'$ of $G$ arises in the following manner. Suppose there is a subgroup $H < G$ such that $\Gbar:=G/H \simeq \Z/ne\Z \times \Z/e\Z$. Then, as we have already noted in the terminal case, $H^2(\Gbar,k^*) \simeq \Z/e\Z$. By Kummer theory, there is a 2-cocycle $\alpha\in H^2(\Gbar,\mu_e)$ which maps onto a generator for $H^2(\Gbar,k^*)$. We then obtain a commutative diagram with exact rows
\[\diagram
1 \rto & \mu_e \dto \rto & G' \dto \rto & G \dto \rto & 1  \\
1 \rto & \mu_e \rto & \Gbar' \rto & \Gbar \rto & 1
\enddiagram\]
where the bottom row comes from the cocycle $\alpha$ and the top row is induced from the bottom by pullback. In brief, the skew group ring $\eps k[[s,t]]*G'$ comes from an element of $H^2(G,k^*)$ which is the image of a generator of $H^2(\Gbar,k^*) \simeq \Z/e\Z$.

The next theorem was proved in [CHI, section~8] and can also be verified using proposition~\ref{pdecomp} and the local computation in proposition~\ref{pskewterminal}. 
\begin{thm}   \label{tskewcanonical}  
The skew group ring $\eps k[[s,t]]*G'$ with $G'$ described in the table below is a canonical order. 
\end{thm}

In the table below, we describe $G$ via listing generators as matrices in $k^{2 \times 2}$. The action of $G$ on $k[[s,t]]$ will be the natural linear action. In all cases except $L_n$, we shall specify the skew group ring $\eps k[[s,t]]*G'$ by writing down the subgroup $H < G$ above. 

\begin{tabular}{c|c|c|c|c}
type  & $G$ & order of $\zeta$ & $H$  & $e$  \\ \hline
$A_{1,2,\z}$ & 
$\begin{pmatrix}
 \z & 0 \\
  0 & 1 
\end{pmatrix} \ \ , \ \ 
\begin{pmatrix}
  1 & 0  \\
  0 & \z
\end{pmatrix}$ & $2e$ & 
$\begin{pmatrix}
  -1 & 0  \\
   0 & -1
\end{pmatrix}$ & $e$  \\\hline
$BL_n$  & 
$\begin{pmatrix}
 \z & 0 \\
  0 & \z^{-1} 
\end{pmatrix} \ \ , \ \ 
\begin{pmatrix}
  0 & 1 \\
  1 & 0
\end{pmatrix}$  & $2n+1$ &
$\begin{pmatrix}
 \z & 0 \\
  0 & \z^{-1} 
\end{pmatrix}$ & 1    \\\hline
$B_n$  & 
$\begin{pmatrix}
 \z & 0 \\
  0 & \z^{-1} 
\end{pmatrix} \ \ , \ \ 
\begin{pmatrix}
  0 & 1 \\
  1 & 0
\end{pmatrix}$  & $2n$ &
$\begin{pmatrix}
 \z & 0 \\
  0 & \z^{-1} 
\end{pmatrix}$ & 1    \\\hline
$L_n$  & 
$\begin{pmatrix}
 \z & 0 \\
  0 & \z^{-1} 
\end{pmatrix} \ \ , \ \ 
\begin{pmatrix}
  0 & 1 \\
  1 & 0
\end{pmatrix}$   & $2n+2$ & see below &   2   \\\hline
$DL_n$  &
$\begin{pmatrix}
 \z & 0 \\
  0 & \z^{-1}
\end{pmatrix} \ \ , 
\ \ 
\begin{pmatrix}
  0 & 1 \\
 -1 & 0
\end{pmatrix} \ , 
\ 
\begin{pmatrix}
 -1 & 0 \\
  0 & 1
\end{pmatrix}$   & $4n-2$ &
$\begin{pmatrix}
 \z & 0 \\
  0 & \z^{-1}
\end{pmatrix} \ , 
\ 
\begin{pmatrix}
  0 & 1 \\
 -1 & 0
\end{pmatrix}$ &  2    \\\hline
$BD_n$  & 
$\begin{pmatrix}
 \z & 0 \\
  0 & \z^{-1}
\end{pmatrix} \ , 
\ 
\begin{pmatrix}
  0 & 1 \\
 -1 & 0
\end{pmatrix} \ , 
\ 
\begin{pmatrix}
 -1 & 0 \\
  0 & 1
\end{pmatrix}$   & $4n-4$ &
$\begin{pmatrix}
 \z & 0 \\
  0 & \z^{-1}
\end{pmatrix} \ , 
\ 
\begin{pmatrix}
  0 & 1 \\
 -1 & 0
\end{pmatrix}$ &  2    \\\hline
$A_{n,\z}$  & 
$\begin{pmatrix}
 \z & 0 \\
  0 & \z^{-1} 
\end{pmatrix} \ \ , \ \ 
\begin{pmatrix}
  1 & 0 \\
  0 & \z^{n+1}
\end{pmatrix}$  & $(n+1)e$ & 
$\begin{pmatrix}
 \z^e & 0 \\
  0 & \z^{-e} 
\end{pmatrix}$  & $e$   \\\hline
$A^1_n,D^1_n,E^1_{6,7,8}$  &  $G < SL_2$  &  & $G$ & 1 \\ \hline
\end{tabular}
\vspace{2mm}

The group $G'$ in the case $L_n$ is the dihedral group of order $8n+8$, that is, 
\[ G'  = \lan \s,\tau | \s^{4n+4} = \tau^2=1, \tau \s = \s^{-1} \tau \ran.\]
Note that $Z(G') = \lan \s^{2n+2}\ran \simeq \mu_2$ and $G'/Z(G')\simeq G$ so $G'$ is indeed a central extension of $G$ by $\mu_2$. Also, the idempotent 
\[ \eps = \frac{1}{2} (1 + \s^{2n+2})   .\]

\section{Skew group rings: minimal resolutions}   \label{sskewres}  

In this section, we recall the definition of a minimal resolution of a canonical order and give conditions for when the recipe in the previous section may be used to construct minimal resolutions as skew group rings. 

Let $Z=\spec R$ where $R$ is a noetherian complete local normal domain and $f:\Zt\lm Z$ be a proper birational morphism. An order $\At$ on $\Zt$ is {\em terminal} if it is terminal complete locally at every closed point and furthermore, it is maximal at the generic point of every exceptional curve. In cases of interest, we assume $f$ to be a rational resolution in the sense that $R^1f_* \calo_{\Zt} = 0$.

Terminal orders on rational resolutions are uniquely determined by their ramification data in the following sense. 

\begin{prop}   \label{ptermunique} 
Let $\Zt \lm \spec R$ be a rational resolution. If $\At,\At'$ are terminal orders on $\Zt$ with the same ramification data, then they are Morita equivalent. 
\end{prop}
\textbf{Proof.} We know from [A87, lemma~1.4] that the Brauer group $ Br \Zt = 0$ is trival so the central simple algebra $\At \otimes_{\Zt} k(\Zt)$ is determined up to Morita equivalence by its ramification data. Replacing $\At,\At'$ by matrix algebras $\At^{n \times n}, (\At')^{m \times m}$ for appropriate $m,n$, we may assume $\At \otimes_{\Zt} k(\Zt) = \At' \otimes_{\Zt} k(\Zt) =: Q$. 

We know that $\At,\At'$ are maximal orders away from a finite set of curves $C_1, \ldots ,C_r$ which are transverse to the exceptional fibre. Also, from definition~\ref{dterminal}, we see that the $C_i$ are mutually disjoint. Consider the localisations $\At_{C_i}, \At'_{C_i}$ at the generic point of $C_i$. We know $\At_{C_i} \simeq \At'_{C_i}$ either by [CI05, proposition~2.11] (or rather, its proof), or by Artin-de Jong's descent theory argument [AdJ, section~1.2]. The Skolem-Noether theorem implies that we can find  invertible elements $q_i \in Q$ such that $\At'_{C_i} = q_i^{-1}\At_{C_i}q_i$. Let $P_i$ be the reflexive hull of $\At q_i \At' \subset Q$ (as an $\calo_{\Zt}$-module). Note that $P_i$ is an $(\At,\At')$-bimodule and that away from all the $C_j$ we have 
\[ \At' = \send_{\At} P_i, \At = \send_{\At'} P_i  \hspace{3cm}  (*) \]   
since $\At,\At'$ are maximal there. Also (*) holds locally at $C_i$ by construction.  We may glue the $P_i$ to obtain a reflexive $\calo_{\Zt}$-module which is equal to $P_i$ locally at $C_i$ and is $P_{i_0}$ elsewhere, for some $i_0$. Since the gluing takes place inside $Q$, it is clear that $P$ is still an $(\At,\At')$-bimodule and 
\[ \At' = \send_{\At} P, \At = \send_{\At'} P   .\]
Finally, $P$ is locally projective as an $\At$-module and an $\At'$-module, so is a Morita bimodule exhibiting the desired Morita equivalence. 

\vspace{2mm}

Consider now a canonical order $A$ with centre $R$. By [CI05, proposition~3.15], $R$ is a log terminal surface singularity and hence, a rational singularity.  Let $f: \Zt \lm \spec R$ be a proper birational morphism which, we note is a rational resolution if $\Zt$ is smooth. We let $f^\# A$ to be a reflexive order on $\Zt$ prescribed at codimension points by the following conditions.

\begin{enumerate}
\item If $C$ is a non-exceptional curve then $(f^\# A)_C = (f^*A)_C$.
\item If $C$ is exceptional then $(f^\# A)_C$ is any maximal order containing $(f^*A)_C$. 
\end{enumerate}

Needless, to say $f^\# A$ is not unique, but the choices corresond precisely to the choices of maximal orders containing $f^* A$ at the various exceptional curves so, in particular, $f^\# A$ exists. Also, all choices have the same ramification data. Indeed, the ramification of $f^\# A$ at any non-exceptional curve is the same as that of $A$, while the ramification at an exceptional curve $E$, is just the ramification of the central simple algebra $A \otimes_R k(R)$ at $E$. The condition that $f^\# A$ is terminal is thus independent of this choice of order. Furthermore, the argument in proposition~\ref{ptermunique} can be used to show that the various choices for $f^\# A$ are all Morita equivalent. 

\begin{defn}  \label{dtermres}  
A {\em terminal resolution} of a canonical order $A$ with centre $R$ is a rational resolution $f:\Zt \lm \spec R$ and a choice of order $\At := f^\# A$ as above which is terminal. We say that such a terminal resolution is {\em minimal} (or is a {\em minimal resolution}) if given any other terminal resolution with centre say $\Zt'$, there exists a morphism $\Zt' \lm \Zt$.  
\end{defn}

Proposition~\ref{ptermunique} shows that minimal resolutions of canonical orders are unique up to Morita equivalence so we shall speak of {\em the} minimal resolution. Also, existence holds by 

\begin{prop} ([CI05, proposition~3.17]) \label{pminexist}  
Any canonical order has a minimal resolution.
\end{prop}

The Artin-Mumford sequence provides an easy recipe for computing the ramification of $f^\# A$ from the ramification of $A$ (see [CK, p.159-160] for an explanation). In particular, it is easy to compute the ramification data of the minimal resolution of a canonical order. This was given in [CHI] and, for the convenience of the reader, the answer will be repeated in section~\ref{scases} when we perform the case-by-case analysis. 

\vspace{2mm}

We now address the question of constructing the minimal resolution of canonical order via skew group rings. To this end, let us consider terminal ramification data $\mathcal{R}=\{\Zt,e_C,e_p\}$ on $\Zt$ where $e_C$ are the primary ramification indices on curves and $e_p$ are the secondary ramification indices at intersection points of ramification curves. In the case of a minimal resolution of a canonical order, we know from [CHI, theorem~6.2] that the exceptional curves occur in one of the following types. 

\vspace{2mm}
\noindent
\textbf{Type $0$}  The exceptional curve $E$ is disjoint from the ramification curves.

\vspace{2mm}
\noindent
\textbf{Type $I_{a,b}$} The ramification index at the exceptional curve $E$ is $e_E = b$. There exist two disjoint ramification curves $U,V$ intersecting $E$ with intersection multiplicity 1 and ramification indices $e_U=e_V=ab$. There are no other ramification curves intersecting $E$. We allow $b=1$ in this case. 

\vspace{2mm}
\noindent
\textbf{Type $C_n$} There is no ramification along the exceptional curve $E$ i.e. $e_E = 1$. There is only one ramification curve $U$ intersecting $E$ and we have $E.U = 2$ and $e_U = n$.

\vspace{2mm}
\noindent
\textbf{Type $X_e$} There is no ramification along the exceptional curve $E$. There are only two ramification curves $U,V$ intersecting $E$. Furthermore, $U$,$V$ intersect each other and $U.E=V.E=1, e_U=e_V=e$.

\vspace{2mm}

The notation was chosen using the fact that the letters $I,C,X$ look like the configuration of ramification curves. 

\begin{prop}  \label{pskewgroup}  
Consider terminal ramification data $\mathcal{R}=\{\Zt,e_C,e_p\}$ on $\Zt$ as above. Suppose there exists a smooth Galois cover $\Xt$ of $\Zt$ with Galois group $\Gbar \simeq \Z/ne\Z \times \Z/e\Z$ such that the ramification indices of $\Xt/\Zt$ are given by the $e_C$. Finally, suppose that any (irreducible) exceptional curve $E \subset \Zt$ is one of the following types. 
\begin{enumerate}
\item $E$ is type $0$.
\item $E$ is type $I_{n,e}$.
\item $E$ is type $C_2$ and $n=2,e=1$.
\item $E$ is type $I_{e,1}$ and $n=1$.
\end{enumerate}
Let $\Gbar'$ be a central extension of $\Gbar$ corresponding to a generator of $H^2(\Gbar,k^*)$. Then $\eps \calo_{\Xt} * \Gbar'$ is a terminal order with ramification data given by $\mathcal{R}$. Here $\eps \in k\mu_e$ is as usual, a primitive idempotent such that a generator $\rho$ for $\mu_e$ acts as a primitive $e$-th root of unity in $\eps k\mu_e$ (see section~\ref{sskew} equation~(\ref{eidemp})).  
\end{prop}
\textbf{Proof.} This follow from the local computations of propositions~\ref{pskewterminal} and \ref{pdecomp}. 

\vspace{2mm}

This will allow us to construct the minimal resolutions of canonical orders via skew group rings in every case except $L_n$ and $DL_n$. Both these cases have exceptional curves of type $X_e$. The more serious obstruction though, is that, they do not have smooth Galois covers whose ramification indices match the ramification of the terminal resolution.

\section{Equivariant derived McKay correspondence} \label{sequi} 

When the minimal resolution of a canonical order is a skew group ring, we may use an equivariant version of the derived McKay correspondence to study it. This correspondence is essentially due to Kawamata and I thank Paul Hacking for bringing my attention to it as well as providing a short proof. We will approach it from the point of view of non-commutative algebra to obtain the version most convenient for us. 

To prove Kawamata's derived equivalence [Kaw], we will use the classical non-equivariant version and so pause a moment to review it. We follow the terminology and setup in [BKR]. Let $H < SL_2$ be a finite subgroup acting naturally on $W:=\spec k[[s,t]]$ and $X:=W/G$ be the corresponding canonical surface singularity. Let $\Xt$ denote the minimal resolution of $X$ which can also be interpreted as the equivariant Hilbert scheme $H-Hilb\ W$ as follows. Recall an {\em $H$-cluster} on $W$ is an $H$-invariant 0-dimensional closed subscheme $T \subset W$ such that the $H$-module $\calo_T \simeq kH$. They are parametrised by the scheme $H-Hilb\ W$. There is a corresponding {\em universal $H$-cluster} $\Wt$ which is a finite cover of $\Xt$. Then Kapranov-Vasserot's derived McKay correspondence [KV, \S 1.4 theorem] (also known to Gonzales-Sprinberg-Verdier) states that $\calo_{\Wt}$ is a tilting bundle inducing a derived equivalence between $\calo_{\Xt}$ and $\ow*H$. By this we mean there is an equivalence of the bounded derived categories of coherent modules $D^b_c(\calo_{\Xt}) \simeq D^b_c(\ow*H)$. 

\begin{prop} (Kawamata) \label{pkawa} 
Let $G'$ be a finite group acting linearly on $W=\spec k[[s,t]]$ via the group homomorphism $\phi:G' \lm GL_2$. Suppose $H$ is a normal subgroup of $G'$ such that $\phi(H) \subset SL_2$. Write $X:=W/H$, $\Xt:= H-\Hilb W$ and $\Wt$ for the universal $H$-cluster. Then $M:=\calo_{\Xt}*(G'/H) \otimes_{\Xt} \calo_{\Wt}$ is a tilting bundle on $\calo_{\Xt}*(G'/H)$ which induces a derived equivalence between  $\calo_{\Xt}*(G'/H)$ and $\End M = \ow*G'$. 
\end{prop}
\textbf{Proof.} We know that $\calo_{\Wt}$ is a tilting bundle on $\Xt$ and that $\ow * H = \End_{\Xt} \calo_{\Wt}$. Note that the $G'$-action on $W$ induces $G'$-actions on $\Wt,\Xt$ which in turn induce a $G'/H$-action on $\Xt$. Now $\calo_{\Xt}*(G'/H)$ is flat over $\calo_{\Xt}$ so adjoint associativity gives 
\begin{equation} 
\Rhom_{\calo_{\Xt}*(G'/H)}(M,N) = \Rhom_{\Xt}(\calo_{\Wt},N) .
\label{eadjoint} 
\end{equation}
If we let $N$ range over arbitrary $\calo_{\Xt}*(G'/H)$-modules and use the fact that $\calo_{\Wt}^{\perp} = 0$ in $D^b_c(\Xt)$, then we deduce that $M^{\perp} = 0$. To show that $M$ is a partial tilting bundle and compute its endomorphism algebra, first note that  
\[ M = \bigoplus_{\bar{g} \in G'/H} \bar{g} \otimes \calo_{\Wt} .\]
If we lift $\bar{g}$ to $g \in G'$,then since $\Wt \lm \Xt$ is $G$-equivariant we have the following isomorphism of left $\calo_{\Xt}$-modules
\[ \bar{g} \otimes \calo_{\Wt} \simeq (\bar{g}^{-1})^* \calo_{\Wt} \simeq (g^{-1})^* \calo_{\Wt} \simeq \calo_{\Wt}.\]
Hence (\ref{eadjoint}) gives the following isomorphism of right modules over $\End_{\Xt} \calo_{\Wt} = \ow * H$.
\begin{equation} \label{eenddec}
 \Rhom (M,M) = \Rhom(\calo_{\Wt},M) \simeq  \bigoplus_{\bar{g} \in G'/H} \ow * H 
\end{equation} 
In particular, we see that $M$ is a tilting bundle. Note that $\ow*G'$ acts on the right of $M$ as follows. First $\ow$ acts on $\calo_{\Wt}$ canonically. Given $g,g' \in G', r \in \calo_{\Wt}$ we define the $G'$-action by  
\[ (g'H \otimes r)g := g'gH \otimes (g^{-1})^* (r) .\]
Note that if $J \subseteq G'$ is a set of coset representatives for $H$, then identifying the summands of $\ow*G' = \oplus_{j \in J} \ow*Hj$ with those in (\ref{eenddec}), we see that $\End M = \ow*G'$ as desired. 

\vspace{4mm}

We will apply the proposition to the groups $G'$ arising in theorem~\ref{tskewcanonical}. Accordingly, we consider the following setup. Let $G\subset GL_2$ act on $W = \spec k[[s,t]]$. Let $H$ be a normal subgroup $G$ which is actually contained in $SL_2$ and write $\Gbar:=G/H$. We consider a central extension $\Gbar'$ of $\Gbar$ by some $\mu_e$. As we have seen, pull-back induces a central extension $G'$ of $G$ by $\mu_e$  we still have $H \simeq \ker (G' \lm \Gbar')$. We may and will thus consider $H$ as a subgroup of $G'$. Now both $\calo_{\Xt}*\Gbar'$ and $\ow*G'$ contain $k\mu_e$ in their centres so any primitive idempotent $\eps$ of $k\mu_e$ gives a central idempotent of $\calo_{\Xt}*\Gbar'$ and $\ow*G'$. 

\begin{cor} \label{ckawa} 
The derived equivalence of the previous proposition induces a derived equivalence between $\e\calo_{\Xt}*\Gbar'$ and $\e\ow*G'$
\end{cor}

%The Artin-Verdier version of the McKay correspondence gives the following lemma which we will use repeatedly. 

%\begin{lemma} \label{lactonexc} \marginpar{lactonexc}
%In the notation of the previous lemma, the action of $g \in G$ on $\Xt$ is given by conjugation by $g$ on $H$.
%\end{lemma} \marginpar{MOVE}
%\textbf{Proof.} Recall that the non-trivial conjugacy classes of $H$ correspond to the non-trivial irreducible representations which in turn correspond to exceptional curves in $\Xt$. We decompose the $H$-module $\ow$ into isotypic components $M_{\nu}$, one for each irreducible representation $\nu \in \hat{H}$. Each component is a (direct sum of) isomoprhic indecomposable Cohen-Macaulay $\ox$-modules. Now $g$ permutes these isotypic components in a manner compatible with the conjugation action. Let $f:\Xt \lm X$ be the natural contraction so that by [AV], the torsion-free part of $f^*M_{\nu}$ is (a direct sum of) $\calo_{\Xt}(C_{\nu})$ where $C_{\nu}$ is the curve transverse to the exceptional curve corresponding to $\nu$. This completes the proof of the lemma. 

\section{Numerical McKay correspondence}   \label{snumMcKay}  

In this section we see how the equivariant derived equivalence of the previous section gives rise to the numerical McKay correspondence observed case by case for canonical orders. Unfortunately, the result does not cover all the canonical orders. 

Our setup for this section will be as follows. Let $G\subset GL_2$ act on $W = \spec k[[s,t]]$. Let $H$ be a normal subgroup of $G$ which is actually contained in $SL_2$. Hence $X:=W/H$ will be a commutative canonical singularity and we let $\Xt$ be its minimal resolution. Suppose that $\Gbar:=G/H \simeq \Z/ne\Z \times \Z/e\Z$ and let $\eps \ow*G'$ be the skew group ring constructed from an extension $G'$ corresponding to a generator of $H^2(\Gbar,k^*)$ as in section~\ref{sskew}. Recall that $H$ can also be treated as a normal subgroup of $G'$ and we set $\Gbar':=G'/H$. As before, we let $\tilde{W}$ be the universal $H$-cluster on $\Xt$. 

One way to extract the classical numerical McKay correspondence from the derived version is to decompose the tilting bundle. We will need this decomposition in the classical case which follows from Artin-Verdier theory [AV] as follows. Below we let $\mathcal{E}$ denote the {\em extended set of exceptional curves} consisting of all the irreducible exceptional curves of $\Xt$ together with the divisor $0$. Also, when we say two curves on $\Xt$ are transverse, we will mean they intersect each other in a point with multiplicity one.

\begin{thm}   \label{tartinverd} ([AV, lemma~1.10])
For each exceptional curve $E\subset \Xt$ we pick a transverse curve $C$ which is disjoint from all the other exceptional curves. We let $F(E)$ denote the universal extension defined by the exact sequence 
\[  0 \lm \Ext^1_{\Xt}(\calo(C),\calo_{\Xt}) \otimes_{\Xt} \calo_{\Xt} \lm F(E) \lm \calo(C) \lm 0 .\]
By default we also set $F(0) = \calo_{\Xt}$. The indecomposable summands of $\calo_{\Wt}$ are precisely the $F(E)$ as $E$ ranges over $\mathcal{E}$. 
\end{thm}
We will use the notation $F(E)$ throughout the rest of the paper. 

Note that $\Gbar$ acts on $\mathcal{E}$ by pull-back of sheaves and each $\Gbar$-orbit of exceptional curves corresponds to an exceptional curve on the centre of $\eps \calo_{\Xt}*G'$. We have the following corollary to proposition~\ref{pkawa}. 
\begin{cor}  \label{cnumMcKay} 
Let $\mathcal{I} \subset \mathcal{E}$ be a set of representatives, one from each $\Gbar$-orbit of $\mathcal{E}$. The indecomposable reflexive $\eps\ow * G'$-modules are the direct summands of 
\[  \Hom_{\eps\calo_{\Xt}*\Gbar'}(\eps\calo_{\Xt}*\Gbar' \otimes_{\Xt} \calo_{\Wt},\eps\calo_{\Xt}*\Gbar' \otimes_{\Xt} F(E))\]
where $E$ ranges over $\mathcal{I}$. 
\end{cor}
\textbf{Proof.} Note that $\eps\ow*G'$ is global dimension two and finite as a module over a complete local noetherian ring. Hence the indecomposable reflexive $\eps\ow * G'$-modules are just the direct summands of $\eps\ow*G'$. Now proposition~\ref{pkawa} shows that these correspond, via the functor \[Hom_{\eps\calo_{\Xt}*\Gbar'}(\eps\calo_{\Xt}*\Gbar' \otimes_{\Xt} \calo_{\Wt},-)\] 
to the direct summands of $\eps\calo_{\Xt}*\Gbar' \otimes_{\Xt} \calo_{\Wt}$. Artin-Verdier's theorem~\ref{tartinverd} shows that $\calo_{\Wt}$ decomposes as the direct sum of the $F(E), E \in \mathcal{E}$. If $E,E'$ lie in the same $\Gbar$-orbit then 
\[  \eps\calo_{\Xt}*\Gbar' \otimes_{\Xt} F(E) \simeq \eps\calo_{\Xt}*\Gbar' \otimes_{\Xt} F(E')  \]
so the corollary follows. 
\vspace{2mm}

For the reader familiar with [CHI], recall that for a canonical order $A$, the number of permissible $A$-modules was one more than the number of exceptional curves in the minimal resolution of $A$ ([CHI, theorem~8.6]). Usually, these permissible modules were just indecomposable reflexive $A$-modules, but occasionally, they were direct sums of a pair of indecomposable reflexive $A$-modules. These two cases correspond to whether or not $\eps\calo_{\Xt}*\Gbar' \otimes_{\Xt} F(E)$ is indecomposable or splits into two indecomposables. The rest of this section is concerned with computing the number of non-isomorphic indecomposable summands of $\eps\calo_{\Xt}*\Gbar' \otimes_{\Xt} F(E)$.

The following lemmas will be useful. 

\begin{lemma}   \label{lFiequiv} 
Suppose that $\Gbar$ stabilises the exceptional curve $E \subset \Xt$ and the image of the induced map $\Gbar \lm \Aut E$ is cyclic. Then there exists a natural $\Gbar$-action on $F(E)$ which makes $F(E)$ an $\calo_{\Xt}*\Gbar$-module.
\end{lemma}
\textbf{Proof.} Since the $\Gbar$-action on $E$ factors through some cyclic quotient group, there are two fixed points on $E$. Looking complete locally at such a fixed point $p \in \Xt$, we can find a curve $C$ which intersects $E$ at $p$ with multiplicity one, such that $C$ is $\Gbar$-stable. This follows from the argument in [A86, proof of proposition~4.13] using the fact that $\Gbar$ stabilises $E$. Note that $\calo(C)$ is a $\Gbar$-stable subsheaf of $k(\Xt)$ so is naturally a $\calo_{\Xt}*\Gbar$-module. 

Suppose first that $E=E_0$ is an extremal curve in the tree of exceptional curves. Since only one other exceptional curve intersects it, we can find a transverse curve $C_0$ to $E_0$ as above, which is $\Gbar$-stable and disjoint from other exceptional curves. Also, $F(E_0)$ is the universal extension of $\calo(C_0)$ by direct sums of $\calo_{\Xt}$ so $F(E_0)$ is also an $\calo_{\Xt}*\Gbar$-module in a natural way. Suppose now that $E$ is not extremal. There are two cases. The first possibility is that $\Gbar$ actually stabilises all the curves in a chain going from $E$ to some extremal exceptional $E_0$ curve in the tree. In this case, we let $C_0$ denote that transverse curve to $E_0$ above and $C$ be a transverse curve to $E$ which intersects no other exceptional curve. Now $C$ is linearly equivalent to a sum of the exceptional curves $E_0,\ldots , E$ and $C_0$. These curves are all $\Gbar$-stable so we can conclude as before. Finally, the only other possibility is that the exceptional fibre is a string of rational curves and $E$ is the middle one. The action of any $g \in \Gbar$ either stabilises all the exceptional curves in the string or reverses the order of the exceptionals. We may assume that some $g\in \Gbar$ swaps the ends of the string. The fixed points of $\Gbar$ on $E$ above therefore cannot be the points of intersection with the neighbouring exceptional curves. Hence, the transverse curve $C$ found in the previous paragraph is disjoint from all the other exceptional curves. The lemma now follows. 

\vspace{2mm}

\begin{lemma}  \label{lbimod} 
Let $F$ be an $\calo_{\Xt}*\Gbar$-module. Then $\calo_{\Xt}*\Gbar' \otimes_{\Xt} F$ is an $\calo_{\Xt}*\Gbar'$-bimodule if one defines the left action to be left multiplication in $\calo_{\Xt}*\Gbar$ and the right action to be given by
\[ (r \otimes_{\Xt} f)g := rg \otimes g^{-1}(f)   \]
where $g \in \Gbar'$, $r$ is a section of $\calo_{\Xt}*\Gbar'$, $f$ is a section of $F$ and the action of $g^{-1}$ on $f$ is given by the $\calo_{\Xt}*\Gbar$-module structure on $F$. Furthermore, $\eps\calo_{\Xt}*\Gbar' \otimes_{\Xt} F$ is an  $\eps\calo_{\Xt}*\Gbar'$-sub-bimodule of $\calo_{\Xt}*\Gbar' \otimes_{\Xt} F$.
\end{lemma}
\textbf{Proof.} To check the above right action makes $\calo_{\Xt}*\Gbar' \otimes_{\Xt} F$ into a bimodule involves checking the relations in the skew group ring. We omit this.  Since $F$ is an $\calo_{\Xt}*\Gbar$-module, $k\mu_e$ acts centrally on the bimodule. This gives the last statement of the lemma. 

\vspace{2mm}

\begin{lemma}  \label{lks} 
Let Add denote the full subcategory of $\eps\calo_{\Xt}*\Gbar'-\mo$, consisting of modules $M$ whose underlying sheaf structure $M_{\calo_{\Xt}}$ is isomorphic to a finite direct sum of $F(E)$ for various $E \in \mathcal{E}$. Then $M \in Add$ is indecomposable if and only if its endomorphism ring is a complete local ring with residue field $k$. In particular, Add is a Krull-Schmidt category, that is, decomposition into indecomposables is unique.
\end{lemma}
\textbf{Proof.} Given $M \in Add$, the derived McKay correspondence for the Kleinian singularity $X$ shows that $End M$ is a finitely generated module over the commutative complete local noetherian ring $\calo_W^G$. The lemma now folows from standard Krull-Schmidt theory over complete local rings (see for example [CR, section~6B]). 

\vspace{2mm}

We seek a numerical form of the McKay correspondence for canonical orders. Since the ramification data of the minimal resolution is the what we know most about the resolution, it will be most useful to find the number of non-isomorphic indecomposable summands of $\eps\calo_{\Xt}*\Gbar' \otimes_{\Xt} F(E)$ in terms of the ramification data of $\eps\calo_{\Xt}*\Gbar'$. Recall that when the skew group ring is chosen as in proposition~\ref{pskewgroup}, the ramification data can also be interpreted as follows. The centre of $\eps \calo_{\Xt} * \Gbar'$ is $\Zt:= \Xt/\Gbar$. Also, if $e_C$ denotes the ramification index of $\Xt/\Zt$ of at $C$, then $e_C$ is also the ramification index of $\eps \calo_{\Xt}* \Gbar'$ at $C$. We assume henceforth that $\Zt$ is smooth, as happens whenever $\eps\calo_{\Xt}*\Gbar'$ is a terminal order. 
 
\begin{thm}   \label{tIcase} 
Let $E\subseteq \Xt$ be an exceptional curve whose image $E_Z \subset \Zt$ is of type $I_{n,e}$. Then there are $n$ non-isomorphic indecomposable summands of $\eps\calo_{\Xt}*\Gbar' \otimes_{\Xt} F(E)$.
\end{thm}
\textbf{Proof.}  Let $\lambda,\nu\in \Gbar'$ be elements which map to generators of $\Gbar$ of order $ne,e$ respectively. As in example~\ref{ecentralext}, we may assume that $\lambda,\nu$ have orders $ne,e$ too. Hence $k\lan \lambda \ran$ has $ne$ primitive idempotents $\eps_j$ and we label them so that $(\lambda - \zeta^j)\eps_j = 0$ where $\zeta$ is a primitive $ne$-th root of unity. Now the ramification of $\Xt/\Zt$ shows that $E\subset \Xt$ is the only exceptional curve lying above $E_Z \subset \Zt$. So by lemma~\ref{lbimod}, $\eps\calo_{\Xt}*\Gbar' \otimes_{\Xt} F(E)$ is an $\eps\calo_{\Xt}*\Gbar'$-bimodule. We thus have a decomposition of left modules 
\[ \eps\calo_{\Xt}*\Gbar' \otimes_{\Xt} F(E) \simeq \bigoplus_{j =0}^{ne-1} M_j \ \ \mbox{where} \ \  
M_j := \eps\calo_{\Xt}*\Gbar' \otimes_{\Xt} F(E) \eps_j   .\]
Now $\nu,\lambda$ skew commute in $\eps k\Gbar'$ by some primitive $e$-th root of unity, so right multiplication by powers of $\nu$ yield isomorphisms
\[ M_j \simeq M_{j+n} \simeq M_{j+2n} \simeq \ldots \simeq M_{j + (e-1)n}   .\]

We wish to show that the $M_j$ are indecomposable, so that the non-isomorphic indecomposable summands of $\eps\calo_{\Xt}*\Gbar' \otimes_{\Xt} F(E)$ are $M_0,\ldots, M_{n-1}$. To this end, note that $\calo_{\Xt}*\lan \lambda \ran$ is a subalgebra of $\eps\calo_{\Xt}*\Gbar'$ so the $M_j$ are $\calo_{\Xt}*\lan \lambda \ran$-modules too. We need the following lemma which gives the $e=1$ case of the proposition. 

\begin{lemma}   \label{lcycliccase} 
Let $F$ be an indecomposable $\calo_{\Xt}*\lan \lambda \ran$-module whose underlying $\calo_{\Xt}$-module structure is a direct sum of copies of $F(E)$. Then 
\[  F \simeq \calo_{\Xt}*\lan \lambda \ran \otimes_{\Xt} F(E) \eps_j  \]
for some $j$. Furthermore, the $\calo_{\Xt}*\lan \lambda \ran \otimes_{\Xt} F(E) \eps_j$ are non-isomorphic $\calo_{\Xt}*\lan \lambda \ran$-modules for $j = 0,1, \ldots , ne-1$. 
\end{lemma}
\textbf{Proof.} Note that the $\calo_{\Xt}*\lan \lambda \ran$-bimodule structure on $\calo_{\Xt}*\lan \lambda \ran \otimes_{\Xt} F(E)$ comes from applying lemma~\ref{lbimod} in the special case $\Gbar = \lan \lambda \ran$. Now the $\calo_{\Xt}*\lan \lambda \ran \otimes_{\Xt} F(E) \eps_j$ are indecomposable, because, as sheaves on $\Xt$, they are isomorphic to $F(E)$. Viewing them as $\lan \lambda \ran$-equivariant sheaves on $\Xt$, we see they differ only in that the action of $\lambda$ differs by an $ne$-th root of unity. Suppose two such modules are isomorphic. Let us fix an $\calo_{\Xt}*\lan \lambda \ran$-module structure on $F(E)$. Then there is some automorphism $\phi$ of the $\calo_{\Xt}$-module $F(E)$ which skew commutes with the action of $\lambda$ by some non-trivial root of unity, that is $\phi \lambda = \xi \lambda \phi$ for some $ne$-th root of unity $\xi\neq 1$. But by lemma~\ref{lks}, $\phi$ is a scalar modulo the radical $J$ of $\End_{\Xt} F(E)$. Now viewing $\lambda$ as an element of $\End_{\Xt/\lan \lambda \ran} F(E)\supset  \End_{\Xt} F(E)$, we see that $\lambda^{-1} \End_{\Xt} F(E) \lambda = \End_{\Xt} F(E)$ and hence $\lambda^{-1} J \lambda = J$. The equation $\phi \lambda = \xi \lambda \phi$ cannot hold in $F(E)/F(E) J$ so the $\calo_{\Xt}*\lan \lambda \ran \otimes_{\Xt} F(E) \eps_j$ must all be non-isomorphic. 

Consider now an indecomposable $F$ as in the lemma. We have a surjective map of $\calo_{\Xt}*\lan \lambda \ran$-modules
\[ \calo_{\Xt}*\lan \lambda \ran  \otimes_{\Xt} F  \lm F\]
which splits as a map of $\calo_{\Xt}$-modules. By Maschke's theorem, it splits as $\calo_{\Xt}*\lan \lambda \ran$-modules. We see that $F$ must be an indecomposable summand of $\calo_{\Xt}*\lan \lambda \ran  \otimes_{\Xt} F(E)$ so the lemma follows from the decomposition
\[ \calo_{\Xt}*\lan \lambda \ran  \otimes_{\Xt} F(E) \simeq \bigoplus_{j=0}^{ne-1} \calo_{\Xt}*\lan \lambda \ran  \otimes_{\Xt} F(E) \eps_j  \]
and the Krull-Schmidt property of lemma~\ref{lks}. 

\vspace{2mm}

We return to the proof of the theorem. We consider the decomposition of each $M_j$ into indecomposable  $\calo_{\Xt}*\lan \lambda \ran$-modules. Note that 
\begin{equation}  \label{exxx}
 M_j = \sum_{i=0}^{e-1} \nu^i \eps \calo_{\Xt}* \lan \lambda \ran \otimes_{\Xt} F(E) \eps_j 
\end{equation}
and in fact, the sum is direct by rank considerations. Now 
\[  \nu^i \eps \calo_{\Xt}* \lan \lambda \ran \otimes_{\Xt} F(E) \eps_j = 
  \eps \calo_{\Xt}* \lan \lambda \ran \otimes_{\Xt} F(E) \nu^i \eps_j  \]
so skew commuting the $\nu^i$ past the $\eps_j$ we see (\ref{exxx}) yields a the following decomposition of $\calo_{\Xt}*\lan \lambda \ran$-modules
\[ M_j \simeq \bigoplus_{i=0}^{e-1} \eps \calo_{\Xt}* \lan \lambda \ran \otimes_{\Xt} F(E) \eps_{j+in}  .\]
This shows that the $M_j$ for $j=0,1, \ldots, n-1$ are already non-isomorphic as $\eps\calo_{\Xt}*\lan \lambda\ran$-modules. It also shows that they are indecomposable as $\eps \calo_{\Xt}*\Gbar'$-modules. This completes the proof of the theorem.

\vspace{2mm}

\begin{prop}   \label{pOCcases}  
Let $E\subseteq \Xt$ be an exceptional curve whose image $E_Z\subset \Zt$ is either i) of type $0$ or ii) of type $C_2$ where $n=2,e=1$. Then $\eps\calo_{\Xt}*\Gbar' \otimes_{\Xt} F(E)$ is indecomposable.
\end{prop}
\textbf{Proof.} In both cases, the $\Gbar$-orbit of $E$ consists of $|\Gbar|$ different exceptional curves. Hence as a sheaf on $\Xt$ we have the decomposition 
\[  \eps\calo_{\Xt}*\Gbar' \otimes_{\Xt} F(E) \simeq \bigoplus_{j} F(E_j)   \]
where the sum runs over the $\Gbar$-orbit of $E$. Suppose now that $F'$ is a non-zero summand of $\eps\calo_{\Xt}*\Gbar'\otimes_{\Xt} F(E)$. The Krull-Schmidt property shows that $F'$ must contain some $F(E_j)$ and hence all $F(E_j)$. It follows that $\eps\calo_{\Xt}*\Gbar' \otimes_{\Xt} F(E)$ is indeed indecomposable.
\vspace{2mm}

\begin{prop}   \label{pXcase} 
Let $E_1\subseteq \Xt$ be an exceptional curve whose image $E_Z\subset \Zt$ is of type $I_{e,1}$ where $n=1$. Then $\eps\calo_{\Xt}*\Gbar' \otimes_{\Xt} F(E_1)$ is the direct sum of $e$ non-isomorphic indecomposable summands. 
\end{prop}
\textbf{Proof.} We note that in this case, $\Gbar = \Z/e\Z \times \Z / e\Z$. The ramification data shows that locally at $E_Z$, the cover $\Xt|_{E_Z}$ of $E_Z$ has $e$ connected components, say $E_1,E_2, \ldots, E_e$, and each component is a cyclic cover of $E_Z$ of degree $e$. Let $\lambda\in \Gbar'$ be such that its image $\bar{\lambda}$ in $\Gbar$ generates the decomposition group of $E_1$. Note that $\Gbar$ is abelian so $\lambda$ stabilises all the $E_i$ and that $\lan \bar{\lambda} \ran \simeq \Z/e\Z$. We may pick $\nu \in \Gbar'$ so that $\lambda,\nu$ generate $\Gbar'$. As in example~\ref{ecentralext}, we may adjust $\lambda, \nu$ by scalars so that $\lambda^e = 1 = \nu^e$. Note also that $\lambda,\nu$ skew commute by some primitive $e$-th root of unity $\zeta$ and $\nu$ permutes the exceptional curves $E_1,E_2, \ldots, E_e$ cyclically. 

Since $\lambda$ stabilises $E_1$, $\eps\calo_{\Xt}*\Gbar' \otimes_{\Xt} F(E_1)$ is an $\eps\calo_{\Xt}*\Gbar'-\calo_{\Xt}*\lan \lambda \ran$-bimodule by lemma~\ref{lbimod}. Hence if $\eps_1, \ldots, \eps_e$ denote the primitive idempotents of $k\lan\lambda\ran$, we have the following decomposition 
\[ \eps\calo_{\Xt}*\Gbar' \otimes_{\Xt} F(E_1) \simeq \bigoplus_{j=1}^e M_j \ \ \mbox{where} \ \ 
M_j:= \eps\calo_{\Xt}*\Gbar' \otimes_{\Xt} F(E_1)\eps_j   .\]
Arguing as in the previous proposition, we see that as a sheaf on $\Xt$ we have
\[  M_j \simeq F(E_1) \oplus F(E_2) \oplus \ldots \oplus F(E_e)   \]
so must be indecomposable. Now lemma~\ref{lcycliccase} shows there are $e$ distinct $\calo_{\Xt}*\lan\lambda\ran$-module structures one can impose on the sheaf $F(E_1)$. They all must occur in $\eps\calo_{\Xt}*\Gbar' \otimes_{\Xt} F(E_1)$ so the $M_j$'s must already be distinct as $\calo_{\Xt}*\lan\lambda\ran$-modules. This finishes the proof of the proposition. 

\vspace{2mm}

The last result we need to determine the indecomposable projectives of $\calo_W * G'$ is the following well-known fact.

\begin{prop}  \label{ptrivialcase} 
$\eps\calo_{\Xt}*\Gbar'$ has $n$ non-isomorphic indecomposable summands. 
\end{prop}
\textbf{Proof.} The proof of theorem~\ref{tIcase} can be repeated, essentially verbatim to give the result in this case. 
\vspace{2mm}

\section{Cyclic covers and minimal resolutions}  \label{sgeneral} 

In this section, we review the construction of orders via non-commutative cyclic covers introduced in [C]. We will construct terminal orders in the complete local case, from this point of view. This will be useful when we construct minimal resolutions of canonical orders. We will also prove that minimal resolutions of canonical orders can always be constructed via non-commutative cyclic covers. A nice feature of the proof is that it is uniform unlike the case-by-case analysis required for the skew group ring construction. 

Below, we will use Artin-Van den Bergh's notion of bimodules (over a scheme), details of which may be found in [AV].

We recall first the definition of non-commutative cyclic covers. Let $Y$ be a normal surface and $\s$ be an automorphism of $Y$ of finite order $e$. Let $G=\lan \s \ran$ and $Z:= Y/G$ which we will assume is a scheme. Consider a rank one reflexive sheaf $L$ on $Y$ which corresponds to some Weil divisor $D$ via $L=\oy(D)$.  We consider the $\oy$-bimodule $L_{\s}$ whose left module structure is $L$, but the right action is twisted through by $\s$ so that $x r = \s(r) x$ for $x \in L, r \in \oy$. One can take tensor products of $\oy$-bimodules and indeed, they form a monoidal category. The only formula we really need involving tensor products will be 

\[ L_{\s}^{\otimes i} = (L \otimes_Y \s^* L \otimes_Y \ldots \otimes_Y L^{\otimes (i-1)})_{\s^i}  .\]

In general, this will not be reflexive as a left (or right) module so we will consider its reflexive hull

\[ L_{\s}^{[i]} := (L \otimes_Y \s^* L \otimes_Y \ldots \otimes_Y L^{\otimes (i-1)})^{**}_{\s^i}  \]
which is also the $\oy$-bimodule $\oy(D + \s^{-1} D + \ldots + \s^{-i+1} D)_{\s^i}$. We call this a {\em reflexive tensor product}. 

Suppose we are given an isomorphism of bimodules $\phi:L^{[e]} \xrightarrow{\sim} \oy$ such that the diagram below commutes.
\begin{equation} \label{eolap}
\begin{CD}
L_{\s} \otimes_Y L_{\s}^{[e-1]} \otimes_Y L_{\s}  @>{1 \otimes \phi}>>
                                                               L_{\s} \otimes_Y \oy \\
@VV{\phi \otimes 1}V                                                  @VV{\phi}V        \\
\oy \otimes_Y L_{\s}            @>{\phi}>>                        L_{\s}
\end{CD}
\end{equation}
In this case, we call $\phi$ a {\em relation satisfying the overlap condition} and we have 
\begin{prop}  \label{pcyclic}  
Suppose we have a relation $\phi:L^{[e]} \xrightarrow{\sim} \oy$ satisfying the overlap condition as above, then there is an algebra structure on 
\[ A(Y;L_{\s},\phi):= \oy \oplus L_{\s} \oplus L_{\s}^{[2]} \oplus \ldots \oplus L_{\s}^{[e-1]}  \]
We call this algebra a {\em non-commutative cyclic cover} of $Y$. It is an order on $Y/G$. 
\end{prop}
\textbf{Proof.} This is proposition~3.2 of [C].  \vspace{2mm}

We shall now construct terminal orders in the complete local case, using non-commutative cyclic covers. Let $\s,\tau$ act linearly on $k[[s,t]]$ via the matrices 
\[ \s = 
\begin{pmatrix}
 \z & 0 \\
  0 & 1 
\end{pmatrix} \ \ , \ \ 
  \tau = 
\begin{pmatrix}
  \z^n & 0 \\
  0 & \z^{-n}
\end{pmatrix} 
\]
where $\z$ is a primitive $ne$-th root of unity. Let $S = k[[s,t]]^{\langle \tau \rangle} = k[[s^e,st,t^e]]$ and note that $\s$ descends to an action on $S$. Let $L$ be the $\z$-eigenspace of $\tau$ acting on $k[[s,t]]$ so that $L = Ss + St^{e-1}$. Note that the indecomposable reflexive $S$-modules are the $\tau$-eigenspaces of $k[[s,t]]$ and that in fact, $L^{[i]}$ is the $\z^i$-eigenspace. The $S$-bimodule $L_{\s}$ can be conveniently interpreted as the submodule $L\s \subseteq k[[s,t]]* \Z/ne$ where the cyclic group $\Z/ne$ is generated by $\s$. 

\begin{prop}   \label{ptermcyclic}  
The natural isomorphism $\phi:L_{\s}^{[e]} \xrightarrow{\sim} S$ is a relation satisyfing the overlap condition. The coresponding non-commutative cyclic cover
\[A:=A(\spec S; L_{\s}) = S \oplus L\s \oplus \ldots \oplus L^{[e-1]}\s^{e-1}  \]
is a terminal order with centre $k[[u,v]]$ where $u=s^{ne},v=t^e$. The ramification curves are $u=0,v=0$ with corresponding ramification indices $ne,e$. 
\end{prop}
\textbf{Proof.}  
The relation $\phi$ is induced by multiplication in $k[[s,t]]*\Z/ne$. Associativity of this mutliplication gives the overlap condition so the cyclic cover is well-defined. Now, we know from [C, corollary~3.4] that $A$ is an order with centre $S^{\lan \s \ran} = k[[s^{ne},t^e]]$. We know from [C, theorem~3.6] that $A$ is a normal order with primary ramification indices equal to the ramification indices of $S/S^{\lan \s \ran}$, that is $ne$ at $u=0$, $e$ at $v=0$ and unramified elsewhere. 

To verify the secondary ramification data, we compute at the ramification curves explicitly. Let $C$ denote the curve $t^e=0$ and $I$ be the ideal $(st,t^e)\triangleleft S$. Note that $I \subset \mbox{rad}\, A_C$. Now
\[ A_C/IA_C \simeq  k(s^e) \oplus k(s^e)s\s \oplus \ldots \oplus k(s^e)s^{e-1}\s^{e-1} \oplus k(s^e)\s^e \oplus k(s^e)s\s^{e+1} \oplus \ldots \oplus k(s^e)s^{e-1}\s^{ne-1} .\]
We see that $A_C/IA_C$ is generated over $k(s^e)$ by $s\s$. One readily computes from this fact that $Z(A_C/IA_C)=k(s^{ne})[(s\s)^n]$. This is a totally ramified field extension of degree $e$ over the residue field $k(C)=k(s^{ne})$. Now the ramification index of $A$ at $C$ is $e$ so $k(s^{ne})[(s\s)^n]$ is also the ramification of $A$ at $C$. Similarly, if $C$ is the ramification curve $s^{ne}=0$ and $I=(s^e,st) \triangleleft S$ we see that 
\[ A_C/IA_C \simeq k(t^e) \oplus k(t^e)t\s^{-1} \oplus \ldots \oplus k(t^e)\s^{-e} \oplus \ldots \oplus k(t^e)t^{e-1}\s^{1-ne}  \]
which is already commutative so is the ramification of $A$ at $C$. This is a degree $ne$ extension of the residue field $k(C)\simeq k(t^e)$ which decomposes as a direct product of $n$ rings corresponding to the $n$ primitive idempotents of $k\lan\s^e\ran \subset A_C/IA_C$. Let $\eps$ be one of these primitive idempotents. Then $\eps A_C/IA_C = \eps k(t^e)[\eps t\s]$ which is a totally ramified field extension of $k(C)$. The preceding two computations show that $A$ is a terminal order with the ramification stated in the proposition. 

\vspace{2mm} 

Once we move away from the complete local case, we find that terminal orders cannot always be constructed via non-commutative cyclic covers. It is quite a surprise that minimal resolutions of canonical orders can always be constructed this way.

\begin{thm}  \label{tcyclicexist} 
Let $A$ be a canonical order on $Z= \spec R$ where $R$ is a complete local ring. Let $\tilde{A}$ be a minimal resolution of $A$ and $\Zt$ its centre. There is a normal surface $\tilde{Y}$ and a non-commutative cyclic cover $\tilde{B}$ of  $\tilde{Y}$ such that $\tilde{B}$ is Morita equivalent to $\tilde{A}$. 
\end{thm}
\textbf{Proof.} We will need to use the canonical bimodule $\w_{\At} := \shom_{\Zt}(\At,\w_{\Zt})$ and the canonical divisor 
\[  K_{\At} := K_{\Zt} + \sum_C (1 - \frac{1}{e_C})C    \]
where the sum ranges over all ramification curves and $e_C$ denotes the ramification index of $\At$ at $C$ as usual. We define $\w_A,K_A$ similarly. More information about $K_A$ can be found in [CI05, \S3.3,3.4] where it plays a central role in the minimal model program for orders. We recall from [CK, proposition~5] that if $\At$ is rank $n^2$ as a sheaf on $\Zt$, then 
\[ \w_{\At}^{[n]} \simeq \At \otimes_{\Zt} \calo_{\Zt}(nK_A) . \]
We will also use the theory of Artin covers developed in [A86, section~2] and further explained in [CHI, section~3]. 

Let $f:\Zt \lm Z$ be the natural contraction which we note is a rational resolution by [CI05, proposition~3.15]. From the original definition of canonical orders [CI05, page~444] we know $K_{\At} \equiv f^*K_A$ so since $K_A$ is $\Q$-Cartier, $K_{\At}$ is numerically trivial. In particular, some multiple of $K_{\At}$ is trivial which in turn implies that some reflexive tensor power of $\omega_{\At}$ is also trivial. (In fact, if $e$ is the index of the central simple algebra $A\otimes_Z k(Z)$ then $\omega_{\At}^{[e]} \simeq \At$.) Similarly, this same reflexive tensor power of $\omega_A$ is also trivial. We may thus form the index one covers $\tilde{B}$ of $\At$ and $B$ of $A$ as defined by Le Bruyn, van Oystaeyen, Van den Bergh (see [CHI, section~4]). In short, they are
\[  \tilde{B} = \At \oplus \w_{\At} \oplus \w_{\At}^{[2]} \oplus \ldots \oplus, \w_{\At}^{[e-1]} \ , 
  \ \ B = A \oplus \w_A \oplus \w_A^{[2]} \oplus \ldots \oplus, \w_A^{[e-1]} . \]
By [CHI, proposition~5.1], $\tilde{B}$ is an Artin cover of $\At$ with the property that its centre $\Yt$ has the same ramification on $\Zt$ as does $\At$. A similar statement can be made for the centre $Y$ of $B$. We also observe from the proof of [CHI, proposition~5.1] that $\Yt,Y$ are regular in codimension one. Furthermore, $\calo_{\Yt},\oy$ are reflexive as sheaves on $\Zt,Z$ respectively. It follows that $\Yt,Y$ are normal by Serre's criterion. 
 
The adjunction formula induces a non-commutative relative trace map $f_* \omega_{\At} \lm \omega_A$ which in turn induces a map $\Yt \lm Y$. We thus obtain the following diagram of surfaces
\[\diagram
 \Yt^{\des} \dto^f \rto & \Yt \dto^f \rto &  \Zt \dto^f  \\
  Y^{\des} \rto       &  Y \rto       &    Z
\enddiagram\]
where the superscript ``res'' denotes the minimal desingularisation which is possible as the surfaces involved are all normal. Now away from the exceptional locus, $\Yt$ and $Y$ are isomorphic so the maps $f$ are all birational.

Now $Y$ is a quotient singularity by [A86, proposition~2.12] and so in particular, is rational. This shows that the closed fibre of $Y^{\des}$ and hence $\Yt^{\des}$ is a tree of rational curves. By [A87, lemma~1.4], $\Yt^{\des}$ has trivial Brauer group. Hence, away from ramification, $\At$ is in the relative Brauer group $\Br \Yt/\Zt$ so by [C, theorem~3.9], its Brauer class is represented by a non-commutative cyclic cover on $Y$. In particular, this cyclic cover has the same ramification data as $\At$ so the result follows from proposition~\ref{ptermunique}. 

\vspace{2mm}

\noindent
\textbf{Remark:} As for the construction of a resolution $\At$ via skew group rings, the key is to find a Galois cover $\Yt$ of $\Zt$ with the same ramification data as $\At$. However, here the cover needs to be cyclic but it need not be smooth.

\vspace{2mm}

\section{Case by case analysis of minimal resolutions}  \label{scases} 

In this section, we construct explicitly, the minimal resolutions of canonical orders via non-commutative cyclic covers and, where possible, via skew group rings. This will be carried out case-by-case, according to ramification data. We then describe the McKay corresondences. In particular, this will prove theorems~\ref{tMcKay} and \ref{tnumMcKay}.

To keep notation straight, we will use the following notation throughout the section. We let $A$ denote a canonical order and $\At$ a minimal resolution. We will have some subset of the following commutative diagram of surfaces, to be described below.
\begin{equation}\diagram
 \Wt \dto^f \rto & \Xt \dto^f \rto & \Yt \dto^f \rto & \Zt \dto^f \\
  W \rto       &  X \rto   & Y \rto     &    Z
\enddiagram
\label{eaux}
\end{equation}
We let $\Zt,Z$ be the centres of $\At,A$ respectively. We let $\Yt$ denote a cyclic cover of $\Zt$ which has the same ramification data as $\At$. In the examples below, there will be a corresponding cyclic cover $Y$ of $Z$. Examples of such surfaces include the surfaces $\Yt,Y$ the proof of theorem~\ref{tcyclicexist}. If $A$ is non-trivial in the Brauer group, then $\Yt$ will not be smooth. Below $\Xt$ will denote (when it exists), a smooth cover of $\Zt$ with the same ramification as $\At$. We let $X$ be the corresponding cover of $Z$. In several examples, $\Xt$ will be constructed from $\Yt$ by taking a Galois cover which is unramified over the smooth locus of $\Yt$. This $\Xt,X$ will coincide with the $\Xt,X$ in sections~\ref{sequi} and \ref{snumMcKay}. In particular, when the smooth cover $\Xt$ exists, we will be able to obtain the derived and numerical McKay correspondences. As in sections~\ref{sequi} and \ref{snumMcKay}, $W = \spec k[[s,t]]$ and $\tilde{W}$ will be the universal $H$-cluster where $H$ is the Galois group of $W/X$. If $G$ denotes the Galois group of $W/Z$, then in all cases we will find that $H \triangleleft G$ and $\Gbar := G/H \simeq \Z/ne \times \Z/e$ for some integers $n,e$. Also, the hypotheses of proposition~\ref{pskewgroup} on exceptional curves will hold, so we may consider the central extensions $\Gbar',G'$ of $\Gbar,G$ by $\mu_e$ as in that proposition. Finally, $\eps$ will denote, as usual, a primitive idempotent of $k\mu_e$ such that a generator $\rho$ of $\mu_e$ acts as a primitive root of unity in $\eps k\mu_e$ (see section~\ref{sskew} equation~(\ref{eidemp})). 

\subsection{Type $A_{1,2,\z}$} \label{stypeA12} 

\textbf{Ramification data of minimal resolution:}
Let $Z = \spec k[[u,v]]$. A type $A_{1,2,\z}$ canonical order has ramification curve the node $uv=0$ with ramification indices $2e$ on the branches and secondary ramification $e$ at the node (see definition~\ref{dcanonical}). Its minimal resolution $\At$ has centre $\Zt$ the blowup of $Z$ at the closed point. The ramification indices are of course $2e$  at the strict transforms of the curves $u=0,v=0$ and is $e$ at the exceptional curve $E_Z\subset \Zt$. This exceptional curve is of type $I_{2,e}$. 

\vspace{2mm}\noindent
\textbf{Construction of minimal resolution:} 
We construct some auxiliary surfaces as in the diagram~(\ref{eaux}) at the beginning of the section.
Let $W = \spec k[[s,t]]$ and $G = \Z/2e \times \Z/2e = \lan\s,\tau\ran$ act on $W$ by $\s:s \mapsto \z s, t\mapsto t, \tau:s \mapsto s, t\mapsto \z t$ where $\z$ is a primitive $2e$-th root of unity. This coincides with the action in the table following theorem~\ref{tskewcanonical}. Note that $Z \simeq W/G$ and $(\s\tau^{-1})^e = \s^e\tau^e$. We define 
\[X :=W/\lan\s^e \tau^e\ran \ \ , \ \ Y = W/\lan \s\tau^{-1}\ran .\] 
Let $\Wt$ be the blowup of $W$ at its closed point and $E$ be its exceptional curve. Note that $G$ also acts on $\Wt$. Hence we obtain a diagram as in (\ref{eaux}) above on setting 
\[\Xt :=\Wt/\lan\s^e \tau^e\ran \ \ , \ \ \Yt = \Wt/\lan \s\tau^{-1}\ran .\]
Note that $\s^e \tau^e$ acts as the antipodal map on $W$ so $X$ is the ordinary double point.
Computing on an affine patch shows that $\s^e \tau^e$ acts on $\Wt$ freely away from $E$ and that $E$ consists of fixed points. The quotient $\Wt/\lan\s^e \tau^e\ran$ is thus smooth with a $(-2)$-curve for an exceptional curve. Hence $\Xt$ is just the minimal resolution of $X$ and $\Wt$ is just the universal $H$-cluster on $W$ where $H=\lan \s^e\tau^e\ran$. Similarly, an affine computation shows $\Zt = \Wt/G$ and $\Wt \lm \Zt$ is ramified precisely along $E$ and the strict transforms of $u=0$ and $v=0$. Moreover, the ramification indices are all $2e$. 

By Galois theory, the quotient $\Gbar:=G/\lan\s^e \tau^e\ran \simeq \Z/2e \times \Z/e$ acts on $\Xt$ with quotient $\Zt$. The ramification indices are $2e$ along the strict transforms of $u=0$ and $v=0$ and $e$ along $E$. This is the same as the ramification of a minimal resolution of a canonical order of type $A_{1,2,\z}$. We are thus in the situation of proposition~\ref{pskewgroup} and can form a skew group ring $\eps \calo_{\Xt}*\Gbar'$ which is the minimal resolution of a type $A_{1,2,\z}$ canonical order. Here $\Gbar'$ is the central extension $\Gbar$ by $\mu_e$ found in proposition~\ref{pskewgroup} and $\eps$ is the idempotent in that proposition. 

We can also construct explicitly the minimal resolution as a non-commutative cyclic cover using $\Yt$. Note that $\Yt/\Zt$ is a cyclic cover with Galois group $G/\lan \s\tau^{-1}\ran$. Note that we may identify this Galois group with $\lan \s \ran$ since the composite $\lan \s \ran \hookrightarrow G \lm G/\lan \s\tau^{-1}\ran$ is an isomorphism. Also, $Y$ is the rational double point of type $A_{2e-1}$. To compute $\Yt$ we consider it as the quotient of $\Xt$ by the cyclic group $\lan \s\tau^{-1} \ran/ \lan\s^e \tau^e\ran \simeq \Z/e$. Working on affine patches we see that $\Yt$ has one exceptional curve and two singularities, each of type $A_{e-1}$. It follows that the minimal desingularisation $\Yt^{\des}$ of $\Yt$ is the minimal desingularisation of $Y$. Let $E_1,\ldots, E_{2e-1}\subset \Yt^{\des}$ be the exceptional curves listed in order. Then $\pi:\Yt^{\des} \lm \Yt$ is just the contraction of all the exceptional curves bar the middle one $E_e$. Now local computations show $\Xt/\Yt$ is unramified in codimension one, so the ramification indices of $\Yt \lm \Zt$ also coincide with that of the minimal resolution of a canonical order of type $A_{1,2,\z}$. Hence, as in theorem~\ref{tcyclicexist}, we can construct the minimal resolution via a non-commutative cyclic cover of $\Yt$. 

To make this explicit, we need only describe the $\oy$-bimodule $L_{\s}$ used to construct the non-commutative cyclic cover $A(\Yt;L_{\s})$. Let $C_1,\ldots, C_{2e-1}$ be curves transverse to $E_i$ and disjoint from $E_j$ whenever $j\neq i$. We define $L:=\calo_{\Yt}(\pi_*C_1 - \pi_*C_{2e-1})$. Now, $\s$ stabilises each of the exceptional curves $E_i$ so $\s^* L \simeq L$ and $L_{\s}^{[2e]} \simeq L^{[2e]}$. We need to show that $L^{[2e]} \simeq \calo_{\Yt}$ so that we can form the non-commutative cyclic cover $A(\Yt;L_{\s})$. Calculating on $\Yt^{\des}$ we find 
\[ 2e(C_1-C_{2e-1}) \sim  (E_1 + 2E_2 + \ldots + (2e-1)E_{2e-1}) - ((2e-1)E_1 + (2e-2)E_2 + \ldots + E_{2e-1}).\]
But $h$ contracts all exceptionals bar $E_e$ so $\pi_* 2e(C_1-C_{2e-1}) \sim 0$. It follows that $L_{\s}^{[2e]} \simeq \calo_{\Yt}$ as desired.

\vspace{2mm}\noindent
\textbf{McKay correspondences:} Proposition~\ref{pkawa} shows that there is a derived equivalence between the minimal resolution $\At:=\eps \calo_{\Xt}*\Gbar'$ and the skew group ring $A:= \eps \calo_W * G'$. Now theorem~\ref{tskewcanonical} shows that $A$ is a canonical order of type $A_{1,2,\z}$ so we have verified theorem~\ref{tMcKay} in the case of type $A_{1,2,\z}$. We may also use the results of section~\ref{snumMcKay} to deduce the numerical form of the McKay correspondence. The unique exceptional curve $E_Z$ in the minimal resolution of $A$ is of type $I_{2,e}$ so, by theorem~\ref{tIcase}, corresponds to precisely two indecomposable reflexive $A$-modules. The regular representation of $\At$ splits in two by proposition~\ref{ptrivialcase}, so also corresponds to two indecomposable reflexive $A$-modules. This explains the observation in [CHI, section~8.1] that there are two permissible modules in this case, each of which splits in two.

\subsection{Type $BL_n$}

\textbf{Ramification data of minimal resolution:} From definition~\ref{dcanonical}, a canonical order $A$ of type $BL_n$ has centre $Z=\spec k[[u,v]]$ and a single ramification curve $C:v^2=u^{2n+1}$ with ramification index 2. Its minimal resolution $\At$ has centre $\Zt$ which minimally resolves the log pair $(Z,C)$. That is, $\Zt$ is the iterated blowup of $Z$ with exceptional fibre, a string of rational curves $E_1, \ldots, E_n$ with $E_n$ a (-1)-curve and all other curves (-2)-curves. Also, $\At$ has a unique ramification curve $D$ of ramification index 2, which is the strict transform of $C$. The iterated blowup $\Zt$ is such that $D$ is smooth and intersects $E_n$ at one point with multiplicity two.

\vspace{2mm}\noindent
\textbf{Construction of minimal resolution:} Let $Y$ be the double cover of $Z$ ramified on $C$ which we note is a Kleinian singularity of type $A_{2n}$. Also, let $\Yt$ be the double cover of $\Zt$ ramified on $D$ which is smooth (and hence we may take $\Xt=\Yt$). Moreover, the inverse image of $E_i$ in $\Yt$ is the union of two smooth rational curves, say $F_i, F_{2n+1-i}$ and we see that $F_1,\ldots F_{2n}$ is a string of (-2)-curves. Note that $Y^{\des}$ has $2n$ exceptional curves too so being minimal we must have $\Yt = Y^{\des}$. Since the canonical order $A$ of type $BL_n$ is trivial in the Brauer group, the minimal resolution can be taken to be $\At = \calo_{\Yt}* \Z/2$ where $\Z/2$ acts via the covering involution $\tau$ of $\Yt/\Zt$. The cyclic cover and skew group ring constructions coincide in this case for we also have $\At = A(\Yt;(\calo_{\Yt})_{\tau})$. 

\vspace{2mm}\noindent
\textbf{McKay correspondences:} To invoke proposition~\ref{pkawa}, we consider the subgroup $G=G'\subset GL_2$ generated by 
\begin{equation} \label{eG}
\s = \begin{pmatrix}
 \z & 0 \\
  0 & \z^{-1} 
\end{pmatrix} \ \ , \ \ 
  \tau = 
\begin{pmatrix}
  0 & 1 \\
  1 & 0
\end{pmatrix} 
\end{equation} 
where $\z$ is a primitive $(2n+1)$-th root of unity. Then $A = \ow*G$ is a canonical order of type $BL_n$ by theorem~\ref{tskewcanonical}. If in proposition~\ref{pkawa} we let $H = \langle \s \rangle$ then we obtain a derived equivalence between $A$ and $\oy * G/H$. Now $G/H \simeq \Z/2$ so to prove theorem~\ref{tMcKay} for type $BL_n$, it suffices to identify the two actions of $\Z/2$ on $\Yt$, one due to the construction of the double cover of $\Zt$ above, the other induced by the action of $\Z/2$ on $Y$. This is clear since both actions are totally ramified on the same curve on $Y$ so coincide on $Y$ and hence on $\Yt$. 

For the numerical McKay correspondence, note that the exceptional curves $E_1,\ldots, E_{n-1}$ are type 0 whilst $E_{n}$ is type $C_2$. Proposition~\ref{pOCcases} shows that they each correspond to an indecomposable reflexive $A$-module. Also, proposition~\ref{ptrivialcase}, shows that the regular representation $\At$ splits in two, so there are two other indecomposable reflexive $A$-modules. This agrees with the observation in [CHI, section~8.2] that there are $n$ indecomposable permissible modules and another permissible module which splits in two. 

\subsection{Type $B_n$}
\textbf{Ramification data of minimal resolution:} This is very similar to the previous case. Again the centre $\Zt$ of the minimal resolution is an iterated blowup of $Z = \spec k[[u,v]]$ so that the exceptional fibre of $f:\Zt \lm Z$ is a string of (-2)-curves $E_1, \ldots ,E_{n-1}$ followed by a (-1)-curve $E_n$ as before. The ramification of the minimal resolution however, consists now of two disjoint smooth curves intersecting $E_n$ transversally. The ramification index of both curves is 2. 

\vspace{2mm}\noindent
\textbf{Construction of minimal resolution:} Let $D$ denote the union of the two ramification curves on $\Zt$. We may form the double cover $\Yt$ of $\Zt$ ramified on $D$ and the double cover $Y$ of $Z$ ramifed on $f(D)$. Note that $Y$ is a Kleinian singularity of type $A_{2n-1}$. Also, $\Yt$ is smooth and the inverse image of $E_i$ in $\Yt$ is the union of smooth rational curves $F_i, F_{2n-i}$ (note that all but $E_n$ split into two rational curves). As for type $BL_n$, $F_1,\ldots F_{2n-1}$ is a string of (-2)-curves and $\Yt = Y^{\des}$. Hence the minimal resolution of $A$ can be taken to be $\At = \calo_{\Yt}* \Z/2$. This coincides with the non-commutative cyclic cover construction. 

\vspace{2mm}\noindent
\textbf{McKay correspondences:}
If we let $G$ be as in the $BL_n$ case (see (\ref{eG})) except that now $\z$ is a primitive $2n$-th root of unity, then $A = \ow * G$ is a canonical order of type $B_n$ by theorem~\ref{tskewcanonical}. Also, setting $H=\lan \s \ran$, we see that $Y = W/H$ and $\At = \calo_{\Yt}*G/H$. Hence, proposition~\ref{pkawa} gives a derived equivalence between $A$ and $\At$, thus proving theorem~\ref{tMcKay} for type $B_n$. 

As for the numerical McKay correspondence, each of the exceptional curves $E_1, \ldots, E_{n-1}$ is of type 0 so correspond to an indecomposable reflexive $A$-module. The exceptional curve $E_n$ is of type $I_{2,1}$ so corresponds to two indecomposable reflexives by theorem~\ref{tIcase}. Finally, the module $\At$ splits in two by proposition~\ref{ptrivialcase} so there are only two other indecomposable reflexives. This agrees with [CHI, section~8.3] where we observed $n-1$ indecomposable permissible $A$-modules and two more permissible $A$-modules, each of which split into two. 

\subsection{Vanishing lemmas}

In the next subsection, we consider canonical orders of type $L_n$. Their minimal resolutions cannot be constructed via skew group rings so the results of section~\ref{sequi} do not apply. Instead, we will need to construct a tilting bundle directly, and this subsection contains some vanishing lemmas to facilitate this. Below, $Y$ denotes a Kleinian singularity and $Y^{\des}$ its minimal resolution. We let $F$ denote the fundamental cycle. 

\begin{lemma} \label{lnoH1} 
Let $L\in \Pic Y^{\des}$. If $H^1(F,L|_F) = 0$ then $H^1(Y^{\des},L) = 0$.  
\end{lemma}
\textbf{Proof.} By Zariski's main theorem and Serre duality, it suffices to show 
\[ H^0(mF,L^*(mF)|_{mF})^* = H^1(mF,L|_{mF})  = 0 . \] 
We do so by induction on $m$. The hypothesis of the lemma gives the $m=1$ case. Consider the exact sequence
\[ 0 \lm L^*_F(F) \lm L^*_{(m+1)F}((m+1)F) \lm L^*_{mF}((m+1)F) \lm 0 \]
where the subscript denotes restriction. By definition of the fundamental cycle, there is an effective divisor $C$, linearly equivalent to $-F$ which is transverse to the exceptional curve. Hence $L^*_{mF}((m+1)F)$ embeds in $L^*_{mF}(mF)$ so has trivial space of global sections. The long exact sequence in cohomology now completes the induction. \vspace{2mm}

We restrict now to the case where $Y$ is a Kleinian singularity of type $A_{2n+1}$ so its minimal resolution is a string of $(-2)$-curves. We let $C_i$ be an irreducible curve transverse to the $i$-th exceptional curve in the string. The subscript $i$ here will be considered modulo $2n+2$. We let $C_0 = 0 = C_{2n+2}$. 

\begin{lemma} \label{lFnoH1} 
Let $i \in [0,n-1], j\in [n+1,2n+2]$ be integers. Then 
\[ H^1(Y^{\des}, \calo(\pm(C_i + C_{n+1} - C_n -C_j))) = 0 .\]
\end{lemma}
\textbf{Proof.} We use lemma~\ref{lnoH1} and the fact that $F \sim -C_1 - C_{2n+1}$. We only prove the $+$ case, the $-$ case being similar. It suffices to show that 
\[ H^0(F, \calo_F(-C_1 - C_{2n+1} - C_i - C_{n+1} + C_n + C_j)) = 0 .\]
Now the pole at $C_n \cap F$ is sandwiched in between the zeros at $C_1 \cap F, C_{n+1} \cap F$ and similarly the pole at $C_j \cap F$ is sandwiched in between the zeros at $C_{n+1} \cap F, C_{2n+1} \cap F$. This gives the desired vanishing. 

\subsection{Type $L_n$}
\textbf{Ramification data of minimal resolution:} The centre $\Zt$ of the minimal resolution is an iterated blowup of $Z = \spec k[[u,v]]$ so that the exceptional fibre of $f:\Zt \lm Z$ is a string of (-2)-curves $E_1, \ldots ,E_{n-1}$ followed by a (-1)-curve $E_n$. There are two ramification curves, which intersect each other transversally and also intersect $E_n$ transversally. The ramification indices are two. 

\vspace{2mm}\noindent
\textbf{Construction of minimal resolution:} 
Let $\Yt$ be the double cover of $\Zt$ ramified on the union $D$ of the ramification curves. Let $Y$ be the corresponding double cover of $Z$ ramified on $f(D)$. As for type $B_n$, $Y$ is an $A_{2n+1}$-singularity. Also the exceptional fibre in $\Yt$ is a string of rational curves  but now $\Yt$ has an ordinary double point at the node in the middle of the string. Let $F_1, \ldots, F_{2n+1}$ be the exceptional curves in $Y^{\des}$ listed in order. Counting exceptional curves and computing their intersection matrix shows that $\Yt$ can be obtained from $Y^{\des}$ by contracting the middle exceptional curve $F_{n+1}$. Let $C_i$ be a curve on $Y^{\des}$ cutting $F_i$ transversally. If $\pi$ denotes the contraction $Y^{\des} \lm \Yt$ then we let $L$ be the reflexive sheaf $\pi_* \calo_{Y^{\des}}(C_{n+1}-C_n)$. Now the covering involution $\tau \in \mbox{Gal}\ Y/Z$ acts on $\Pic Y^{\des}$ by swapping $C_i,C_{2n+2-i}$. Furthermore, $(1+\tau)(C_{n+1}-C_n) = -C_n + 2C_{n+1} - C_{n+2} \sim F_{n+1}$ so $L_{\tau}$ is 2-torsion. Hence, we may form the non-commutative cyclic cover $\At = A(\Yt;L_{\tau}) = \calo_{\Yt} \oplus L_{\tau}$. Locally at the ordinary double point on $\Yt$, this construction is the same as that in proposition~\ref{ptermcyclic}, so $\At$ is indeed the minimal resolution of a type $L_n$ canonical order.

\vspace{2mm}\noindent
\textbf{McKay correspondences:}
We cannot apply proposition~\ref{pkawa} in this case. Let $G'=\lan \s,\tau\ran$ be the group defined in theorem~\ref{tskewcanonical} for type $L_n$. Note that $G'$ acts on $W = \spec k[[s,t]]$ via the matrices in  (\ref{eG}) where $\z$ is a primitive $(2n+2)$-th root of unity. Also theorem~\ref{tskewcanonical} shows that $A:= \e \ow*G'$ is a canonical order of type $L_n$ where $\eps = \frac{1}{2}(1 + \s^{2n+2})$. The relationship with the minimal resolution stems from the fact that $Y \simeq W/\lan \s \ran$ and the action of $G'/\lan \s \ran \simeq \lan \tau \ran$ is the Galois action on $Y/Z$. 

\begin{thm}\label{tmckayf}
There is a tilting bundle for $\At$ whose endomorphism algebra is the canonical order $\e \ow * G'$. 
\end{thm}
\textbf{Proof.} For simplicity, we shall let $C_i$ also denote the corresponding Weil divisor on $\Yt$ and the subscript $i$ is considered modulo $2n+2$. Similarly, we set $C_0 = 0 (=C_{2n+2})$. Let $\Wt$ be the universal $\lan \s \ran$-cluster on $W$ so that $\calo_{\Wt}$ is a tilting bundle for $Y^{\des}$. 

For $j \in [n+2,2n+2]$ we let  $M_j:= \calo_{\Yt}(C_j)$ and for $i \in [1,n+1]$ we let $M_i: = \calo_{\Yt}(C_{i-1} + C_{n+1} -C_n)$. The choice of notation is made so that the left module $(L_{\tau} \otimes_{\Yt} M_r)^{**} \simeq M_{2n+3 - r}$ and 
\[ M_i \simeq \calo_{\Yt}(C_i + E_i + E_{i+1} + \ldots + E_n) .\]
Finally we let $M = \oplus_{r=0}^{2n+1} M_r$ and wish to show $(\At \otimes_{\Yt} M)^{**}$ is a tilting bundle for $\At$. 

We first show $\Ext^i _{\At}((\At\otimes_{\Yt} M)^{**},(\At\otimes_{\Yt} M)^{**}) = 0$ for $i>0$. Now $\At$ is terminal and $(\At\otimes_{\Yt} M)^{**}$ is a reflexive sheaf, so it is also locally projective as an $\At$-module. Hence, from the local-global Ext spectral sequence we need to show that $H^i(\Yt,\shom_{\At}((\At\otimes_{\Yt} M)^{**},(\At\otimes_{\Yt} M)^{**})) = 0$ for $i>0$. Now $\At\otimes_{\Yt} M \lm (\At\otimes_{\Yt} M)^{**}$ has finite length cokernel so functoriality of the reflexive hull and adjoint associativity show
\begin{equation}\label{ef1ncase}
\shom_{\At}((\At\otimes_{\Yt} M)^{**},(\At\otimes_{\Yt} M)^{**}) = \shom_{\Yt}(M, (\At\otimes_{\Yt} M)^{**}) = \shom_{\Yt}(M, M^{\oplus 2})
\end{equation} 
Now the last term decomposes as a sum of terms of the form $\shom_{\Yt}(M_r,M_s)$ which is a reflexive rank one sheaf which can be written as $\pi_{*} N_{rs}$ for some line bundle $N_{rs}$ on $Y^{\des}$ with $R^1\pi_* N_{rs} = 0$. By the Leray spectral sequence, it remains to show that all the $H^1(Y^{\des},N_{rs}) = 0$ which follows from the fact that $\calo_{\Wt} = \oplus \calo_{Y^{\des}}(C_i)$ is a tilting bundle and from lemma~\ref{lFnoH1}. 

Our next goal is to compute $\End (\At\otimes_{\Yt} M)^{**}$. From (\ref{ef1ncase}) we see $\End (\At\otimes_{\Yt} M)^{**} = (\End_{\Yt} M)^{\oplus 2}$ so it seems expedient to determine $\End_{\Yt} M$.  Note $\pi_*\calo_{\Wt} \subset M$ and generically, they are isomorphic. In fact, the right $\pi_*\calo_{\Wt}$-module structure on $\pi_*\calo_{\Wt}$ extends to $M$ as follows. For $i,l \in [1,2n+1]$ such that $i+l < 2n+2$ we have the following linear equivalences on $Y^{\des}$.  
\[ \begin{split}
C_i & \sim iC_1 + (i-1)E_1 + (i-2)E_2+ \ldots + E_{i-1}  \\
C_l & \sim lC_1 + (l-1)E_1 + (l-2)E_2+ \ldots + E_{l-1}  \\
C_{i+l} & \sim (i+l)C_1 + (i+l-1)E_1 + (i+l-2)E_2+ \ldots + E_{i+l-1}  \\
\end{split}\]
Hence $C_{i+l} - C_i - C_l \sim a_1 E_1 + \ldots + a_{i+l-1}E_{i+l-1}$ for some positive integers $a_1, \ldots, a_{i+l-1}$. This linear equivalence induces the multiplication map $\calo(C_i) \otimes_{\Yt} \calo(C_l) \lm \calo(C_{i+l})$ in $\calo_{\Wt}$. The fact that $a_l,\ldots,a_{i+l-1} \geq 1$ imply that this multiplication extends to a map $\calo(C_i) \otimes_{\Yt} M_l \lm M_{i+l}$. This gives the $\pi_*\calo_{\Wt}$-module structure on $M$. We seek to show that $\End M = \End \calo_{\Wt}$. To this end, note
\[ \Hom_{\Yt}(M_r,M_s) = H^0(Y^{\des}, \calo(C_s - C_r + E_{rs}))  \]
where $E_{rs}$ is some exceptional divisor. Now the classical McKay correspondence on $Y$ shows $\End \calo_{\Wt} = \ow * \lan \s \ran$ so $H^0(Y^{\des}, \calo(C_s - C_r)) = H^0(Y^{\des},\calo(C_{s-r}))$. We claim 
\[H^0(Y^{\des}, \calo(C_s - C_r + E_{rs})) = H^0(Y^{\des}, \calo(C_s - C_r))  \]
so that $\Hom_{\Yt}(M_r,M_s)$ is just given by right multiplication by $H^0(Y^{\des},\calo(C_{s-r}))\subset \ow$. Note that $E_{rs}$ is either an effective divisor or negative effective. In the first case, we have equality because the right hand side is a reflexive $\oy$-module and so determined on the punctured spectrum which is isomorphic to $Y^{\des} - F$. In the second case, it follows since we have checked that $M$ is stable under right multiplication by $\pi_*\calo_{\Wt}$. It is clear now that $\End M = \ow * \lan \s \ran$. 

Now $\End (\At \otimes_{\Yt} M)^{**}$ is generated by $\End M$ and the involution $\tau$ induced by the natural isomorphism $M \simeq (L_{\tau} \otimes_{\Yt} M)^{**}$. To identify this with the canonical order $\e \ow * G'$, it is convenient to rescale the natural action of $\s$ so that it acts by $\xi^{2i-1}$ on $M_i$ where $\xi$ is a square root of $\zeta$. Note that now $\s^{2n+2} = -1$ which is compatible with the algebra structure in $\eps \ow * G'$ since $\eps = \frac{1}{2}(1 + \s^{2n+2})$. Also $(L_{\tau} \otimes_{\Yt} M_r)^{**} \simeq M_{2n+3 - r}$  so $\tau $ swaps $M_r$ and $M_{2n+3 - r}$. This shows $\s\tau = \tau \s^{-1}$. Finally, $\tau$ commutes with $\ow$ in the same fashion as in $\eps \ow * G'$ so we are done. 

We show now that $(\At \otimes_{\Yt} M)^{**}$ does indeed generate the derived category. For convenience write $\At(C)$ for $(\At \otimes_{\Yt} \calo(C))^{**}$. We construct some $\At$-modules in the triangulated subcategory $\mathcal{T}$ generated by the $\At(C_j)$ for $j\in[n+2,2n+2]$. Let $N$ be an invertible sheaf on $\Yt$ and $r \neq n+1$. Note that the exact sequences of sheaves on $\Yt$
\[ 0 \lm N(-C_r) \lm N \lm \calo_{C_r} \lm 0 ,\]
remain exact when tensored with $\At$. Hence if $C$ is any linear combination of $C_j$'s then $\At(-C)\in \mathcal{T}$. 

Consider also the exact sequence of $\At$-modules
\[  0 \lm \At \lm \At(C_{n+1}) \lm Q \lm 0 ,\]
with $Q$ chosen appropriately. It remains exact on tensoring on the right by an invertible $\calo_{\Yt}$-module $N$. Also, $Q$ is supported on the affine subscheme $C_{n+1}$ so we get an exact sequenece
\[  0 \lm \At\otimes_{\Yt} N \lm \At(C_{n+1})\otimes_{\Yt} N \lm Q \lm 0 .\]
Now $\mathcal{T}$ contains $\At(C_i + C_{n+1} - C_n)$ and $\At(C_{n+1} - C_n)$ so it must contain $\At(-C)$ for $C$ any linear combination of $C_1, \ldots C_n,C_{n+2}, \ldots C_{2n+1}$. 
Hence, it suffices to show that for any non-zero $\At$-module $F$, there is a non-zero homomorphism from such a $\At(-C)$ to $F$. This holds by ampleness of $C_1 + \ldots + C_n + C_{n+2} + \ldots + C_{2n+1} $ in $\Yt$. This completes the proof of the theorem.

\vspace{2mm}

\begin{cor}
There are $n+1$ indecomposable reflexive $A$-modules. 
\end{cor}
\textbf{Proof.} This is also contained in [CHI, section~8.4], and we wish to give a proof here via the derived McKay correspondence above. As usual, the indecomposable reflexives correspond to the non-isomorphic indecomposable summands of the tilting bundle $(\At \otimes_{\Yt} M)^{**}$. But these are precisely the $\At \otimes_{\Yt} \calo_{\Yt}(C_j)$ for $j \in [n+2,2n+2]$. There are $n+1$ of these. 

\vspace{2mm}

Note that all exceptional curves are of type 0 or $X_2$ so one easily verifies theorem~\ref{tnumMcKay} in this case.

\subsection{Type $BD_n$}
\textbf{Ramification data of minimal resolution:} As usual, the centre $\Zt$ of the minimal resolution is an iterated blowup of $Z = \spec k[[u,v]]$ so that the exceptional fibre of $f:\Zt \lm Z$ is a string of (-2)-curves $E_1, \ldots ,E_{n-1}$ followed by a (-1)-curve $E_n$. The ramification curves of the minimal resolution this time,  are $E_1, \ldots, E_{n-1}$ together with three smooth curves intersecting $E_1,E_{n-1},E_n$ transversally. The ramification index of all curves is 2. 

Note that the exceptional curves are of type $0$ or $X_2$ so theorem~\ref{tnumMcKay}  holds by inspection. 

\vspace{2mm}\noindent
\textbf{Construction of minimal resolution:} 
Let $D\subset \Zt$ be the union of the ramification curves. Let $\Yt$ be the double cover of $\Zt$ ramified along $D$ which exists by the proof of theorem~\ref{tcyclicexist} (its existence can also be checked easily by direct computation). Let $Y$ be the corresponding double cover of $Z$ ramified on $f(D)$. The exceptional fibre of $\Yt$ is a string of $n$ rational curves but there are ordinary double points in $\Yt$ sitting above all the nodes in the discriminant $D$. Now $Y$ is a Kleinian singularity of type $D_{2n}$ so counting exceptional curves shows that $Y^{\des}$ is also the minimal resolution of $\Yt$. We denote the exceptional fibres in $Y^{\des}$ by $F_1, \ldots , F_{2n}$ so that $F_{2i}$ is the strict transform of $E_i$ and the other exceptionals arise from resolving the ordinary double points in $\Yt$. Note that $F_1, \ldots , F_{2n-2}$ is a string of (-2)-curves and that $F_{2n-1},F_{2n}$ both intersect $F_{2n-2}$ transversally to form the ``fork'' of the $D_{2n}$ Dynkin diagram. We let $C_i$ denote transverse curves to $F_i$ and let $\pi:Y^{\des} \lm \Yt$ be the natural contraction. We seek to form a non-commutative cyclic cover using the reflexive $\calo_{\Yt}$-module 
\[ L := \pi_* \calo_{Y^{\des}}(C_1 - C_2 + C_3 - \ldots + C_{2n-1}) .\]
First note that the covering involution $\rho \in \mbox{Gal}\ Y/Z$ acts trivially on $\Pic Y^{\des}$. We need a relation to form the cyclic cover which will follow if we can show $L^{[2]} \simeq \calo_{\Yt}$. This amounts to showing in $\Pic Y^{\des}$ that 
\[2(C_1 - C_2 + C_3 - \ldots + C_{2n-1}) \in \Z F_1 + \Z F_3 + \ldots + \Z F_{2n-1}.\]
Now computing intersection numbers gives 
\[ F_1 \sim -2C_1+C_2, F_3 \sim C_2-2C_3+C_4, \ldots, F_{2n-1} \sim C_{2n-2}-2C_{2n-1} \]
so 
\[ 2(C_1 - C_2 + C_3 - \ldots + C_{2n-1}) \sim -(F_1 + F_3 + \ldots + F_{2n-1}).\]
We observe thus that $L^{[2]} \simeq \calo_{\Yt}$ and so can form a cyclic cover $\At = A(\Yt;L_{\rho}) = \calo_{\Yt} \oplus L_{\rho}$. A local inspection using proposition~\ref{ptermcyclic} shows that the ramification is indeed that of the minimal resolution of a canonical order of type $BD_n$. 

Note that $L$ can also be used to construct a smooth double cover $\Xt$ of $\Yt$ and hence a minimal resolution can be constructed via a skew group ring construction. The cover $\Xt/\Yt$ is \'etale away from the double points so it is easy to compute the number of exceptional curves of $\Xt$. The adjunction formula also gives their self-intersection numbers from which we deduce that $\Xt$ is the minimal resolution of a $D_{n+1}$ singularity. Note that as in [A, lemma~4.2], $\Xt/\Zt$ is Galois and has Galois group $\Gbar:=\Z/2 \times \Z/2$. Picking a central extension $\Gbar'$ of $\Gbar$ by $\mu_2$ and idempotent $\eps \in k\mu_2$ as in proposition~\ref{pskewgroup}, we see that the minimal resolution can also be constructed as $\At = \eps \calo_{\Xt}*\Gbar'$. 

\vspace{2mm}\noindent
\textbf{McKay correspondences:} To obtain the McKay correspondence using proposition~\ref{pkawa}, we will need to interpret the Galois cover $\Xt/\Zt$ in a different fashion. To this end, let $G$ be the group generated by 
\begin{equation} \label{ebigG}
\s = 
\begin{pmatrix}
 \z & 0 \\
  0 & \z^{-1}
\end{pmatrix} \ \ , 
\ \ \tau = 
\begin{pmatrix}
  0 & 1 \\
 -1 & 0
\end{pmatrix} \ \ , 
\ \ \rho = 
\begin{pmatrix}
 -1 & 0 \\
  0 & 1
\end{pmatrix} 
\end{equation}
where $\z$ is a primitive $(4n-4)$-th root of unity. Note that $Y$ is the quotient of $W$ by $\lan \s,\tau \ran$. We let $X:= W/H$ where $H:=\lan \s^2,\tau \ran$ so that $X$ is a type $D_{n+1}$ canonical singularity.
Note that its minimal resolution is isomorphic to $\Xt$. Furthermore, the action of $G/H\simeq\Z/2\times \Z/2$ on $\Xt$ coincides with the action of the Galois group of $\Xt/\Zt$ since, they agree on $X$ by inspection of the ramification data. It follows that the minimal resolution $\At = \eps \calo_{\Xt}*\Gbar'$ is derived equivalent to $\eps \ow*G'$ which is a canonical order of type $BD_n$ by theorem~\ref{tskewcanonical}. Theorem~\ref{tMcKay} is thus verified in this case. 

%We start by proving $\Yt \simeq \Xt/\lan \s \ran$. First note that $\s$ acts on the conjugacy classes of $H$ by permuting the ones containing $\tau$ and $\tau^{-1}$ while leaving the others invariant. By lemma~\ref{lactonexc}, $\s$ acts on the exceptional curves in $\Xt$ by permuting the two at the fork and leaving the others invariant. Also, $\s$ must fix all the nodes in the exceptional curve but is free away from the exceptional curves as it acts freely on $X$ away from the singular point. Hence the fixed points of $\s$ are either nodes or exceptional curves. The singularities of $\Xt/\lan\s\ran$ are thus ordinary double points. The exceptional curves on $\Xt/\lan \s \ran$ corresponding to (pointwise) fixed lines have self-intersection $-4$ while the ones containing an isolated fixed point will have self-intersection -1. Resolving singularities we see that the minimal resolution of $\Xt/\lan \s \ran$ must be $Y^{\des}$ so $\Xt/\lan \s \ran$ has no curves of self-intersection $<-2$. In particular, $\s$ cannot have any fixed lines. We may thus work out which curves are contracted in $Y^{\des} \lm \Xt/\lan \s \ran$ and hence conclude that $\Yt \simeq \Xt/\lan \s \ran$. 

%We need to show that $\s$ actually maps to the non-trivial element of $\mbox{Gal}\ \Xt/\Yt$. (CHECK this OK).
%it suffices now to show that, up to automorphisms of $\Yt$, there is a unique morphism $\Xt \lm \Yt$. 

%The argument in the type $BL_n$ case, shows that $\rho$ acts on $\Xt$ in a way compatible with the action of $\mbox{Gal}, \Yt/\Zt$. Hence the two actions of $G/H$ coincide. (FIND BETTER PROOF!!)

We also have the numerical McKay correspondence. The exceptional curves $E_1,\ldots,E_{n-1}$ are of type $I_{1,2}$ so each corresponds to a unique indecomposable reflexive $A$-module by theorem~\ref{tIcase}. On the other hand $E_n$ is type $I_{2,1}$ so corresponds to two indecomposable reflexives. Finally, the module $\At$ gives one more indecomposable reflexive module by proposition~\ref{ptrivialcase}. This agrees with the $n$ indecomposable permissible modules found in [CHI, section~8.6] and the additional permissible module which splits in two. 

\subsection{Type $DL_n$}

\textbf{Ramification data of minimal resolution:} As usual, the centre $\Zt$ of the minimal resolution is an iterated blowup of $Z = \spec k[[u,v]]$ so that the exceptional fibre of $f:\Zt \lm Z$ is a string of (-2)-curves $E_1, \ldots ,E_{n-1}$ followed by a (-1)-curve $E_n$. The ramification curves of the minimal resolution this time,  are $E_1, \ldots, E_{n-1}$ together with two curves, one intersecting $E_1$ transversally and the other intersecting $E_{n-1}$ and $E_n$ transversally at their point of intersection. The ramification indices of all curves are 2.

\vspace{2mm}\noindent
\textbf{Construction of minimal resolution:} Let $D$ be the union of ramification curves. We form the double cover $\Yt$ of $\Zt$ ramified on $D$ as in theorem~\ref{tcyclicexist} and let $Y$ be the corresponding double cover of $Z$ ramified on $f(D)$. Sitting above each of $E_1,\ldots ,E_{n-1}$ is an exceptional curve of $\Yt$. Sitting above $E_n$ are two exceptional curves intersecting transversally. As in the previous case, there are ordinary double points in $\Yt$ sitting above all the nodes in the discriminant $D$. Now $Y$ is a Kleinian singularity of type $D_{2n+1}$ so counting exceptional curves shows that $Y^{\des}$ is also the minimal resolution of $\Yt$. We denote the exceptional fibres in $Y^{\des}$ by $F_1, \ldots , F_{2n+1}$ so that for $i \in [1,n-1],\ F_{2i}$ is the strict transform of $E_i$ and $F_1, F_3, \ldots ,F_{2n-1}$ are  exceptionals arising from the ordinary double points in $\Yt$. The exceptional curves $F_{2n},F_{2n+1}$ come from the inverse image of $E_n$. Note that $F_1, \ldots , F_{2n-2}$ is a string of (-2)-curves and that $F_{2n+1},F_{2n}$ both intersect $F_{2n-1}$ transversally. We let $C_i$ denote transverse curves to $F_i$ and let $\pi:Y^{\des} \lm \Yt$ be the natural contraction. The covering involution $\rho \in \mbox{Gal}\ Y/Z$ now acts non-trivially. It fixes $F_1, \ldots, F_{2n-1}$ but switches $F_{2n}$ and $F_{2n+1}$. We seek to construct a cyclic cover using 
\[ L:= \pi_* \calo_{Y^{\des}}(C_1 - C_2 + C_3 - \ldots - C_{2n}) .\]
This time we need to show that $(L \otimes \rho^* L)^{**} \simeq \calo_{\Yt}$ or in other words 
\[(1 + \rho)(C_1 - C_2 + C_3 - \ldots - C_{2n}) \in \Z F_1 + \Z F_3 + \ldots + \Z F_{2n-1}.\]
A calculation similar to the $BD_n$ case shows that 
\[ (1+\rho)(C_1 - C_2 + C_3 - \ldots - C_{2n}) \sim -(F_1 + F_3 + \ldots + F_{2n-1})\]
so again we may form the cyclic algebra $\At = A(\Yt:L_{\rho}) = \calo_{\Yt} \oplus L_{\rho}$. It gives the correct answer by the local computation of proposition~\ref{ptermcyclic}. 

\vspace{2mm}\noindent
\textbf{McKay correspondences:} We do not know if there is a derived equivalence between the minimal resolution and the canonical order of type $DL_n$ in theorem~\ref{tskewcanonical}. However, [CHI, section~8.5] shows that there are $n+1$ indecomposable permissible modules in this case and no others. The exceptional curves $E_1,\ldots, E_{n-1}$ are of type $I_{1,2}$ whilst $E_n$ is of type $X_2$. Hence theorem~\ref{tnumMcKay} is verified in this instance too.  

\subsection{Type $A_{n,\z}$}   \label{stypeAnz}  

\textbf{Ramification data of minimal resolution:} The centre $\Zt$ of the minimal resolution is the minimal resolution of a Kleinian singularity of type $A_n$. Let $E_1, \ldots, E_n$ be the exceptional curves written in order. The ramification curves are $E_1,\ldots, E_n$ together with two other curves, one transverse to $E_1$ and the other transverse to $E_n$. The ramification indices are all $e$ this time where $e$ is the order of $\z$. 

\vspace{2mm}\noindent
\textbf{Construction of minimal resolution:} Let $D$ be the union of ramification curves and $\Yt$ be the $e$-fold cover of $\Zt$, totally ramified on $D$. Also, let $Y$ be the corresponding $e$-fold cover of $Z$ ramified on $f(D)$ which we note is an $A_{ne+e-1}$-singularity. Also the exceptional fibre in $\Yt$ is a string of $n$ rational curves and now $\Yt$ has $A_{e-1}$-singularities at the nodes of the exceptional fibre. Again we have $Y^{\des}$ is the minimal resolution of $\Yt$ (on counting and noting $n+(n+1)(e-1) = ne+e-1$). The exceptional curves in $Y^{\des}$ form a string $F_1, \ldots, F_{ne+e-1}$ where $F_{ie}$ is the strict transform of $E_i$. Let $C_i$ be a curve on $Y^{\des}$ cutting $E$ transversally and $\pi$ denote the contraction $Y^{\des} \lm \Yt$. Again the Galois group of $Y/Z$ acts trivially on $Y^{\des}$ so we seek an $e$-torsion Weil divisor in $\mbox{Cl}\ \Yt$. Let 
\[ L:= \pi_* \calo_{Y^{\des}}(\sum_{i=0}^n C_{ie+1} - \sum_{i=1}^n C_{ie}) .\]
Let $J$ denote the subgroup of $\Pic Y^{\des}$ generated by $\{F_j| e \nmid j\}$. To show $L^{[e]} \simeq \calo_{\Yt}$ it suffices to verify that $e(\sum_{i=0}^n C_{ie+1} - \sum_{i=1}^n C_{ie}) \in J$. We compute 
\begin{align*}
(e-1)F_1 + (e-2)F_2 + \ldots + F_{e-1} & \sim -eC_1 + C_e  \\
(e-1)F_{ie+1} + (e-2)F_{ie+2} + \ldots + F_{ie+e-1} & \sim (e-1)C_{ie}-eC_{ie+1} + C_{ie+e}  \\
(e-1)F_{ne+1} + (e-2)F_{ne+2} + \ldots + F_{ne+e-1} & \sim (e-1)C_{ne}-eC_{ne+1}  
\end{align*}
for $i \in [1,n-1]$. Summing the above equations for all $i$ shows that $L$ is indeed $e$-torsion and that the corresponding cyclic algebra $A(\Yt;L_{\s}) = \calo_{\Yt} \oplus L_{\s} \oplus \ldots \oplus L_{\s}^{[e-1]}$ is well-defined. It has the correct ramification locally so, by proposition~\ref{ptermcyclic} gives the minimal resolution. As for types $BD_n$, we can also use $L$ to form a smooth $e$-fold cover $\Xt$ of $\Yt$ and so construct the minimal resolution $\At$ as a skew group ring $\eps \calo_{\Xt}*\Gbar'$ as in proposition~\ref{pskewgroup}. Note that the exceptional curves of $\Xt$ have self-intersection -2 and so $\Xt$ is isomorphic to the minimal resolution of a type $A_n$ singularity. As in the type $BD_n$ case, $\Xt/\Zt$ is Galois with Galois group $\Z/e \times \Z/e$. 

\vspace{2mm}\noindent
\textbf{McKay correspondences:} 
We use proposition~\ref{pkawa} with $G \subset GL_2$ the group generated by
\[ \s = 
\begin{pmatrix}
 \z & 0 \\
  0 & \z^{-1} 
\end{pmatrix} \ \ , \ \ 
  \tau = 
\begin{pmatrix}
  1 & 0 \\
  0 & \z^{n+1}
\end{pmatrix} 
\]
where $\z$ is a primitive $(n+1)e$-th root of unity. We let $H= \lan \s^e \ran$ so that $\Gbar:=G/H \simeq \Z/e \times \Z/e$. Note that the minimal resolution of $X:=W/H$ is indeed $\Xt$. We lift the central extension $\Gbar'$ of $\Gbar$ to a central extension $G'$ of $G$ as in section~\ref{sskew}. Note that by theorem~\ref{tskewcanonical}, $A:=\eps \ow * G'$ is a canonical order of type $A_{n,\z}$. Furthermore, proposition~\ref{pkawa} shows it is derived equivalent to $\At=\eps \calo_{\Xt} * \Gbar'$. 

All the exceptional curves here are of type $I_{1,e}$ so theorem~\ref{tIcase} and \ref{ptrivialcase} show that there are $n+1$ indecomposable reflexive $A$-modules. This matches with the $n+1$ indecomposable permissible modules found in [CHI, section~8.8]. 

%As in the $BD_n$ case, we wish to identify the above action of $\Gbar$ with the action given in proposition~\ref{pkawa}. I don't have a good proof of this fact FIX!! We consider the action given in proposition~\ref{pkawa} first. Now $G$ is abelian so lemma~\ref{lactonexc} shows that $\s$ leaves each of the exceptional curves in $\Xt$ invariant. All the nodes are fixed points of $\s$ and the only other fixed points that can occur are if there are whole exceptional curves which are fixed pointwise. I think this should force $\Xt/\lan \s \ran \simeq \Yt$ and we can conclude equality of actions in the type $BD_n$ case. 

\vspace{1cm}
\textbf{\large{References}}

\begin{itemize}
\item {[{A86}]} M. Artin, ``Maximal orders of global dimension and 
                   Krull dimension two'', Invent. Math., \textbf{84}, (1986), p.195-222
\item{[{A87}]}  M. Artin, ``Two-dimensional orders of finite representation type''  Manuscripta Math.  \textbf{58},  (1987),  p.445--471.
\item {[{AdJ}]} M. Artin. J. de Jong, ``Stable Orders over Surfaces'',
 in preparation.
\item {[{AVer}]} M. Artin, Verdier, ``Reflexive modules over rational double points'', Math. Ann. \textbf{270}, (1985), p.79-82
\item  {[{AV}]} M. Artin, M. Van den Bergh,  ``Twisted Homogeneous Coordinate Rings'',
        J. of Algebra, \textbf{133}, (1990) , p.249-271
\item {[{BKR}]} T. Bridgeland, A. King, M. Reid, ``The McKay correspondence as an equivalence of derived categories'', Journal of the AMS \textbf{14} (2001), p.535-554
\item {[{C}]} D. Chan,``Noncommutative cyclic covers and maximal orders on surfaces'', Advances in Math.  \textbf{198} (2005), p. 654-83
\item {[{CHI}]} D. Chan,P. Hacking, C. Ingalls,``Canonical singularities of orders over surfaces'' submitted, available on arXiv:math/0401425
\item {[{CI04}]} D. Chan, C. Ingalls,``Noncommutative coordinate rings and stacks'', Proc. of  LMS vol. 88 (2004), p.63-88
\item {[{CI05}]} D. Chan, C. Ingalls, ``The minimal model program for orders over surfaces'', Invent. Math. \textbf{161} (2005), p.427-452
\item {[{CR}]}   R. Curtis, I. Reiner, ``Methods of representation theory. Vol. I. With applications to finite groups and orders'', John Wiley \& Sons,  New York, (1981)
\item {[{Kaw}]}  Y. Kawamata, ``Equivalences of derived categories of sheaves on smooth stacks'', Amer. J. of Math. \textbf{126} (2004), p. 1057-83
\item {[{KV}]} E. Vasserot, M. Kapranov, ``Kleinian singularities, derived categories and Hall algebras'',  Math. Ann.  \textbf{316}  (2000), p.565--576
\item {[{V}]} Van den Bergh, ``Three dimensional flops and noncommutative rings'', Duke Math. J. \textbf{122}  (2004), p.423--455
\end{itemize}

\end{document}